\documentclass[12pt]{article}
\usepackage[a4paper]{geometry}

\usepackage{amsthm} 
\usepackage{amsmath}
\usepackage{amssymb}
\usepackage{amsfonts}
\def\N{\mathbb{N}}

\def\P{\mathbb{P}}

\def\R{\mathbb{R}}
\def\E{\mathbb{E}}
\def\L{\mathbb{L}}

\def\Z{\mathbb{Z}}

\def\S{\mathcal{S}}

\def\D{\operatorname{\mathcal{D}iv}}

\def\J{\mathcal{J}}
\def\G{\mathcal{G}}
\def\<{\big\langle}
\def\>{\big\rangle}

\def\sym{{\operatorname{sym}}}

\def\Vol{{\operatorname{Vol}}}
\def\Div{{\operatorname{div}}}

\def\eref#1{(\ref{#1})}

\newtheorem{Theorem}{Theorem}[section]

\newtheorem{Lemma}[Theorem]{Lemma}

\newtheorem{Proposition}[Theorem]{Proposition}

\newtheorem{Condition}[Theorem]{Condition}
\theoremstyle{remark}
\newtheorem{Remark}[Theorem]{Remark}
\theoremstyle{definition}

\newtheorem{Definition}[Theorem]{Definition}

\begin{document}

\title{Approximation of the effective conductivity of ergodic media by periodization.\protect\footnotetext{AMS 1991 {\it{Subject Classification}}. Primary 74Q20, 37A15; secondary   37A25} \protect\footnotetext{{\it{Key words and phrases}}. Effective conductivity, periodization of ergodic media, Weyl decomposition}}         
\author{Houman Owhadi\footnote{
LATP, UMR CNRS 6632, CMI, Universit\'{e} de Provence , owhadi@cmi.univ-mrs.fr }}        
\date{\today}          
\maketitle
\abstract{
This paper is concerned with the approximation of the effective conductivity $\sigma(A,\mu)$ associated to an elliptic operator $\nabla_x A(x,\eta) \nabla_x$ where for $x\in \R^d$, $d\geq 1$, $A(x,\eta)$ is a bounded elliptic random symmetric $d\times d$ matrix and $\eta$ takes value in an ergodic probability space $(X,\mu)$. Writing $A^N(x,\eta)$ the periodization of $A(x,\eta)$ on the torus $T^d_N$ of dimension $d$ and side $N$ we prove that for $\mu$-almost all $\eta$
$$
\lim_{N\rightarrow +\infty}\sigma(A^N,\eta)=\sigma(A,\mu)
$$
We extend this result to non-symmetric operators $\nabla_x (a+E(x,\eta)) \nabla_x$ corresponding to diffusions in ergodic divergence free flows ($a$ is $d\times d$ elliptic symmetric matrix and $E(x,\eta)$ an ergodic skew-symmetric matrix); and to discrete operators corresponding to random walks on $\Z^d$ with ergodic jump rates.\\
The core of our result is to show that the ergodic Weyl decomposition associated to $\L^2(X,\mu)$ can almost surely be approximated by  periodic Weyl decompositions with increasing periods, implying that  semi-continuous variational formulae associated to $\L^2(X,\mu)$ can almost surely be approximated by variational formulae minimizing on periodic potential and solenoidal functions.
} 

\tableofcontents

\section{Introduction}
Homogenization theory has been developed to find the asymptotic behavior of operators associated to an heterogeneous ergodic medium when the microscopic scale associated to the heterogeneities tends towards $0$ in front of the macroscopic scale of the observation. The mathematical formulation of this theory \cite{BeLiPa78} has been   first developed in the simpler case of elliptic and parabolic periodic operators. The first rigorous results on elliptic and stationary parabolic ergodic operators were obtained by S. Kozlov \cite{Ko80}, \cite{Ko85}, G. Papanicolaou and S. Varadhan \cite{PaVa79} in the late seventies. Next C. Kipnis and S.R.S. Varadhan \cite{KiVa86} followed by \cite{MaFe89} and \cite{OsSa95} introduced powerful  central limit theorems allowing the extension of homogenization theory to a wide range of ergodic operators.\\
Thus two main categories of problems have been addressed by homogenization theory: the asymptotic behavior of  periodic  operators and the asymptotic behavior of  ergodic  operators. The question of the existence of a natural and continuous link between those two categories of applications has naturally arisen. Indeed for large deviations \cite{DGI00} and equilibrium fluctuations \cite{GOS01} of $\nabla \phi$ interface models it has been observed that the regularity of the effective conductivity associated to the infinite dimensional ergodic system under its finite dimensional periodic approximations hides   
an hard core difficulty in extending the mathematical description of relaxation towards equilibrium of periodic environments to ergodic ones.\\
Recently this regularity property  has been proven for the self-diffusion coefficient for the exclusion process  \cite{LOV00}. The case of the effective diffusivity of a symmetric random walk on $\Z^d$, under the condition that its  jump rates are i.i.d. has been addressed in \cite{CapIof01}, which also put into evidence an  exponential rate of convergence of effective diffusivities of the finite volume approximations of the ergodic medium.\\
It is important to note that the periodic approximation of the effective conductivity of a continuous ergodic operator can be obtained from almost sure $G$-convergence properties associated to the homogenization of that operator; this is the so called "principle of periodic localization" (equation 5.15, page 155 of \cite{JiKoOl91}, \cite{Piat02}).\\
It is well known that these effective conductivities associated to periodic and ergodic operators are given by variational formulae in the geometrical framework of Weyl decomposition. Weyl decomposition plays a central role in homogenization in periodic or ergodic media but is independent of them. Thus it is natural to seek for a geometrical property inherent to Weyl decomposition allowing the almost sure approximation of ergodic effective conductivities by periodization.
The purpose of our paper is to show that such property does exist: the Weyl decomposition is gifted with almost sure strong stability properties (theorem \ref{jhdchdcbcc221}, proposition \ref{dddhdhdhhh8711}) which can be used to establish a natural and continuous link between homogenization in periodic and ergodic media. We refer to theorem \ref{JHlkhbsljbj81} for a continuous operator, theorem \ref{dkhbdlhl8981} for a non symmetric operator and \ref{dkhbdlhl8981nn2} for discrete operator (one can also consider a larger class of homogenization problems such as those associated to the lemma 3.1 of \cite{Nor97}). This property inherent to Weyl decomposition can also be used to obtain almost sure periodic approximation results for a wider class of variational formulations on an ergodic space than those associated to an effective conductivity, this is the object of theorem \ref{sddhdvdhdzhvdzd78712}.

\section{General set up}
\subsection{The ergodic space}
\subsubsection{Continuous case}
Let $(X,\G,\mu)$ be a probability space with $\eta \in X$ labeling the particular realization of the quenched medium. We assume that on $(X,\G,\mu)$ acts ergodically a group of measure preserving transformations $G=\{\tau_x\,:\, x\in \R^d\}$, i.e. that the following are satisfied:
\begin{Condition}\label{wekjhbkjb891}
 $\forall x \in \R^d$, $\tau_x$ preserves the measure, namely, $\forall A \in \G$, $\mu(\tau_x A)=\mu(A)$;
\end{Condition}
\begin{Condition}\label{wekjhbkjb892}
 The action of $G$ is ergodic, namely, if $A=\tau_x A$ $\forall x \in \R^d$, then $\mu(A)=0$ or $\mu(A)=1$.
\end{Condition}
Let $L^2(\mu)$ be the Hilbert space of square integrable functions on X with the usual scalar product
\begin{equation}
\int_{X} f(\eta)g(\eta)\,d\mu(\eta)
\end{equation}
Let $f\in L^2(\mu)$, for almost every $\eta$ we define
\begin{equation}
(T_x f)(\eta)=f(\tau_{-x}\eta)
\end{equation}
We assume furthermore that 
\begin{Condition}\label{wekjhbkjb893}
For any measurable function $f(\eta)$ on $X$, the function $T_x f(\eta)$ defined on the Cartesian product $X\times \R^d$ is also measurable (where $\R^d$ is endowed with the Lebesgue measure).
\end{Condition}
It follows that  that $T_x$ form a strongly continuous unitary group on $L^2(\mu)$ (see \cite{JiKoOl91} chapter 7).\\
For $f\in L^1(\mu)$ we write. 
\begin{equation}
<f>\equiv \int_X f(\eta)\mu(d\eta)
\end{equation}

\subsubsection{Discrete case}
We shall distinguish through this paper two cases of ergodic spaces. The one mentioned above associated to a continuous of measure preserving transformations and the one mentioned here associated to a discrete measure preserving transformations.
We shall keep the same notation used above for the continuous case. $(X,\G,\mu)$ will remain our ergodic probability space with $\eta \in X$ labeling the particular realization of the quenched medium but we replace the group of measure preserving transformations acting ergodically on $(X,\G,\mu)$ by $G=\{\tau_x\,:\, x\in \Z^d\}$. We will replace the conditions \ref{wekjhbkjb891} and \ref{wekjhbkjb892} by
\begin{Condition}\label{wekjhbkjb8901}
 $\forall x \in \Z^d$, $\tau_x$ preserves the measure, namely, $\forall A \in \G$, $\mu(\tau_x A)=\mu(A)$;
\end{Condition}
\begin{Condition}\label{wekjhbkjb8902}
 The action of $G$ is ergodic, namely, if $A=\tau_x A$ $\forall x \in \Z^d$, then $\mu(A)=0$ or $\mu(A)=1$.
\end{Condition}

\subsection{Weyl Decomposition}\label{jdhdbbhdhii101}
\subsubsection{Continuous case}
A vector field $f=(f_1,\ldots,f_d)$, $f_i \in L^2_{loc}(\R^d)$, $i=1,\ldots,d$ is called vortex-free in $\R^d$ if
\begin{equation}\label{ksdjhhvd1}
\int_{\R^d}\big(f_i \partial_j \phi-f_j \partial_i \phi\big)dx=0\,\quad \forall \phi \in C_0^\infty(\R^d)
\end{equation}
It is well known that any vortex-free vector possesses a potential function, i.e., admits the representation $f=\nabla u$, $u\in H^1_{loc}(\R^d)$. Therefore  the potentiality of a vector field $f$ is equivalent to the property \eref{ksdjhhvd1}. A vector field $f$ is said to be solenoidal in $\R^d$ if
\begin{equation}
\int_{\R^d}f_i \partial_i \phi(x)\,dx=0,\quad \forall \phi \in C_0^\infty(\R^d)
\end{equation}
Now let us consider vector fields on $X$. A vector field $f\in (L^2(\mu))^d=\L^2(X,\mu)$ will be called potential (resp., solenoidal), if almost all its realizations $T_x f(\eta)$ are potential (resp., solenoidal) in $\R^d$. The spaces of potential and solenoidal vector fields denoted by $\L^2_{pot}(X,\mu)$ and $\L^2_{sol}(X,\mu)$, form closed sets in $\L^2(X,\mu)$.\\
Set
\begin{equation}
F^2_{pot}=\{f\in \L^2_{pot}(X,\mu),\;<f>=0\}
\end{equation}
\begin{equation}
F^2_{sol}=\{f\in \L^2_{sol}(X,\mu),\;<f>=0\}
\end{equation}
By Weyl's decomposition (see the lemma 7.3 of \cite{JiKoOl91}) the following orthogonal decomposition are valid
\begin{equation}\label{aslkjhbiqsuu71}
\L^2(X,\mu)=F^2_{pot}\oplus F^2_{sol}\oplus \R^d=F^2_{pot}+\L^2_{sol}(X,\mu)
\end{equation}

\subsubsection{Discrete case}
For any $f:\, X \rightarrow \R$ and  $i\in \{1,\ldots,d\}$ we write
\begin{equation}
D_i f(\eta)=f(\tau_{-e_i}\eta)-f(\eta)\quad \quad D_i^* f(\eta)=D_i f(\tau_{e_i} \eta)
\end{equation}
Write 
\begin{equation}
\L^2(X,\mu):=\{(f_i)_{1\leq i\leq d}\,:\, f_i \in L^2(\mu)\}
\end{equation}
and $F^2_{pot}$ the completion of $\{(D_i f(\eta))_{1\leq i\leq d}\,:\, f \in L^2(\mu)\}$ in $\L^2(X,\mu)$ with respect to the standard $L^2$ norm ($\|f\|^2=\sum_{i=1}^d\<f_i^2\>$).\\
Write $\S(X,\mu)$ the set of skew-symmetric matrices $H$ such that $H_{i,j} \in L^2(\mu)$ and define $\D H$ as the vector $(\D H)_i:=\sum_{j=1}^d D_j H_{ij}$. We write $F^2_{sol}$ the completion of $\{\D H\,:\, H\in  \S(X,\mu)\}$ 
in $\L^2(X,\mu)$ with respect to the standard $L^2$ norm.\\
We will prove the subsection \ref{skasjsbjbd981}  the following theorem corresponding to the Weyl decomposition.
\begin{Theorem}\label{suvsiuoswzvzvd81}
One has
\begin{equation}\label{suvsiuoswzvzvd80001}
\L^2(X,\mu)=F^2_{pot}\oplus F^2_{sol} \oplus \R^d
\end{equation}
\end{Theorem}

\section{Main result}
\subsection{Almost sure continuity of Weyl's decomposition}
The results and notations given in this section are valid in the continuous case as in the discrete case. Write $T^d_N$ the torus of side $N$ and dimension $d$ ($T^d_N:=\R^d/(N \Z^d)$ in the continuous case and $T^d_N:=\Z^d/(N \Z^d)$ in the discrete).
  Write $\L^2(T^d_N)$ the space of square integrable vectors fields on $T^d_N$ and gift it with the norm $\|f\|_{\L^2(T^d_N)}^2$ (written below) 
to obtain an Hilbert space. We define
$$\|f\|_{\L^2(T^d_N)}^2:=\sum_{i=1}^d N^{-d}\int_{T^d_N}(f_i(x))^2\,dx$$
in the continuous case and
$$\|f\|_{\L^2(T^d_N)}^2:=\sum_{i=1}^d N^{-d}\sum_{x\in T^d_N}(f_i(x))^2$$
in the discrete case.
 In the continuous case, write $F^2_{pot}(T^d_N)$ ($F^2_{sol}(T^d_N)$) the completion of the space of smooth $T^d_N$-periodic potential (solenoidal) forms (with $0$-Lebesgue mean value) in $\L^2(T^d_N)$ with respect to that norm. In the discrete case we shall use the following definitions
\begin{equation}
F^2_{pot}(T^d_N):=\{\nabla f\,:\, f\in L^2(T^d_N)\}
\end{equation}
Where $L^2(T^d_N)$ is the space of square integrable functions on $T^d_N$ and $\nabla$ is the discrete gradient on $\Z^d$, $(\nabla f)_i:=(\nabla_i f)= f(x+e_i)-f(x)$.
For $N$, let us write $\S(T^d_N)$  the set of  skew-symmetric matrices with  coefficients in  $L^2(T^d_N)$ and for $H\in \S(T^d_N)$, $\Div H$ is the vector field defined by $(\Div H)_i= \sum_{j=1}^d \nabla_j H_{i,j}$.\\
\begin{equation}
F^2_{sol}(T^d_N):=\{\Div H\,:\, H\in \S(T^d_N)\}
\end{equation}
As in the ergodic case it is easy (see \cite{JiKoOl91}) to obtain   the following Weyl decomposition for the periodic case:
\begin{equation}\label{aslsssdkjhbiqsuaau71}
\L^2(T^d_N)=F^2_{pot}(T^d_N)\oplus F^2_{sol}(T^d_N)\oplus \R^d
\end{equation}
For any $\xi \in \L^2(X,\mu)$ we write $\Pi_N \xi\in \L^2(T^d_N)$ its periodization on the torus $T^d_N$: for $x= p N+y$ with $p\in \Z^d$ and $y\in [0,N(^d$ ( $y\in \Z\cap [0,N(^d$ in the discrete case)
$$\Pi_N \xi(x,\eta):= T_{y}f(\eta)$$
Let us define $\lim_{N\rightarrow \infty }F^2_{pot}(T^d_N)$ as the subset of $\xi\in\L^2(X,\mu)$ such that 
for $\mu$-almost all $\eta\in X$ there exists a sequence $\nu^N_{pot} \in F^2_{pot}(T^d_N)$ such that
\begin{equation}
\lim_{N\rightarrow \infty} \|\Pi_N \xi(\eta)- \nu^N_{pot}\|_{\L^2(T^d_N)}=0
\end{equation}
Similarly we define $\lim_{N\rightarrow \infty }F^2_{sol}(T^d_N)$ as the subset of $\xi\in\L^2(X,\mu)$ such that 
for $\mu$-almost all $\eta\in X$ there exists a sequence $\nu^N_{sol} \in F^2_{sol}(T^d_N)$ such that
\begin{equation}
\lim_{N\rightarrow \infty} \|\Pi_N \xi(\eta)- \nu^N_{sol}\|_{\L^2(T^d_N)}=0
\end{equation}
The following theorem (valid for both continuous and discrete cases) is the central point linking homogenization in ergodic media to homogenization in periodic media.
\begin{Theorem}\label{jhdchdcbcc221}
\begin{equation}\label{djhfgfgf8817222}
\lim_{N\rightarrow \infty }F^2_{pot}(T^d_N)=F^2_{pot}
\end{equation}
and
\begin{equation}\label{ddddddhhjfhfhfh8811}
\lim_{N\rightarrow \infty }F^2_{sol}(T^d_N)=F^2_{sol}
\end{equation}
Where the limits in $N$ are taken along $\R^+$ in the continuous case and along $\N$ in the discrete.
\end{Theorem}
For any $\xi \in \L^2(X,\mu)$ we write $\xi_{pot}$ and $\xi_{sol}$ its components on $F^2_{pot}$ and $F^2_{sol}$. For any $\nu \in \L^2(T^d_N)$ we write $\nu_{pot}$ and $\nu_{sol}$ its components on $F^2_{pot}(T^d_N)$ and $F^2_{sol}(T^d_N)$. 
The theorem \ref{jhdchdcbcc221} is based on the following proposition
\begin{Proposition}\label{dddhdhdhhh8711}
 For any $\xi \in \L^2(X,\mu)$, for $\mu$-almost all $\eta\in X$ 
\begin{equation}\label{dkdjdjdhjhd090334}
\lim_{N\rightarrow \infty} \|\Pi_N \big(\xi_{pot}(\eta)\big)- \big(\Pi_N \xi(\eta)\big)_{pot}\|_{\L^2(T^d_N)}=0
\end{equation}
and
\begin{equation}\label{dkdjdjdhjhd09033445}
\lim_{N\rightarrow \infty} \|\Pi_N \big(\xi_{sol}(\eta)\big)- \big(\Pi_N \xi(\eta)\big)_{sol}\|_{\L^2(T^d_N)}=0
\end{equation}
Where the limits in $N$ are taken along $\R^+$ in the continuous case and along $\N$ in the discrete.
\end{Proposition}
In order to prove \ref{dkdjdjdhjhd09033445} observe that it is sufficient to prove the following lemma
\begin{Lemma}\label{dddhdhdhhh8711b}
for any $\xi \in \L^2(X,\mu)$, for $\mu$-almost all $\eta\in X$ 
\begin{equation}\label{dkdjdjdhjhd090334b}
\lim_{N\rightarrow \infty} \|\Pi_N \big(\xi_{pot}(\eta)\big)- \big(\Pi_N \xi_{pot}(\eta)\big)_{pot}\|_{\L^2(T^d_N)}=0
\end{equation}
and
\begin{equation}\label{dkdjdjdhjhd09033445b}
\lim_{N\rightarrow \infty} \|\Pi_N \big(\xi_{sol}(\eta)\big)- \big(\Pi_N \xi_{sol}(\eta)\big)_{sol}\|_{\L^2(T^d_N)}=0
\end{equation}
Where the limits in $N$ are taken along $\R^+$ in the continuous case and along $\N$ in the discrete.
\end{Lemma}
There exists two strategies two prove lemma \ref{dddhdhdhhh8711b}.  For the equation \eref{dkdjdjdhjhd090334b} for instance, the second strategy consists in taking the primitive of $\xi_{pot}$, then approximating that primitive by a $T^d_N$-periodic function whose gradient is used to approximate $\Pi_N \big(\xi_{pot}(\eta)\big)$. This second strategy is the  harder one because the primitive of $\xi_{pot}$ is in general not an element of $\L^2(X,\mu)$.  The first strategy consists in finding a function $u$ of $L^2(\mu)$ such that its gradient is an approximation of $\xi_{pot}$ then proving \eref{dkdjdjdhjhd090334b} for $\nabla u$. Our proof of the Weyl decomposition (theorem \ref{suvsiuoswzvzvd81}) and the definition of $F^2_{pot}$ allows us to use the simpler first strategy in the discrete case. In the continuous case we have used the second strategy (which can also be used for the discrete case) but one can also adapt the second strategy to the continuous case (it is important to note that first, one has to prove that there exists a subset of $\L^2(X,\mu)$ such that the gradients of its elements are dense in $F^2_{pot}$).
The proof of lemma \ref{dddhdhdhhh8711b} will be given in subsection \ref{sdkjuhdjhdhh89112}. We will first prove in subsection \ref{skasjsbjbd984} equation \eref{dkdjdjdhjhd09033445b} in the discrete case using the first strategy, the proof of equation \eref{dkdjdjdhjhd090334b} in the discrete case being similar we will just give the idea of the proof in subsection \ref{skasjsbjbd983}. Then in subsection \ref{eicnboiecbb81} we will prove equation \eref{dkdjdjdhjhd090334b} in the continuous case using the second strategy (the proof of equation \eref{dkdjdjdhjhd09033445b} being similar we will just give its idea in subsection \ref{aslslssuxuc1}).
\begin{Remark}
Let us note that the proof of lemma \ref{dddhdhdhhh8711b} with the second strategy is constructive and  Birkhoff ergodic theorem is applied only to $\xi_{pot}$ and not  to a sequence of its approximations. Thus the second strategy associate the rate at which  Birkhoff ergodic theorem holds for $\xi_{pot}$ to the rate at which the limits in lemma \ref{dddhdhdhhh8711b} hold.
More precisely, writing for $M\in \N$, $I(M)=\{(i_1,\ldots,i_d)\in \{1,\ldots, M\}^d\,;\, \min_{j}\min(i_j-1,M-i_j)=0 \}$ and $\J(M)$  the set of cubes $B_i$ indexed  by $i\in I(M)$ and
\begin{equation}\label{sjbsbdcjbdfufb8811sssaa12}
B_i=\{x\in [0,1]^d\,:\, \max_j |x_j-(i_j-0.5)/M|\leq 1/(2M)\}
\end{equation}
one can define
\begin{equation}\label{sjbsssaasbdcjbdfufb8811sssaa12}
\begin{split}
f(N,\xi_{pot},\eta)=\sup \Big\{&M\in \N\,:\, \sup_{B_i\in \J(M^2)} M (\Vol(B_i))^{-1}\big|\int_{B_i} \xi_{pot}(Nx,\eta)\,dx\big|\leq 1\quad\text{and}\\& \sup_{B_i\in \J(M)} (\Vol(B_i))^{-1}\int_{B_i} \big|\xi_{pot}(Nx,\eta)\big|^2\,dx\leq 2\big<\xi_{pot}^2\big>\Big\}
\end{split}
\end{equation}
It is then easy to check by Birkhoff ergodic theorem \ref{wlkkiudbu1} that a.s. $f(N,\xi_{pot},\eta)\rightarrow \infty$ as $N\rightarrow \infty$
and by taking $M=f(N,\xi_{pot},\eta)$ and $P=M$ in subsection \ref{eicnboiecbb81} one obtains that for all $\xi_{pot}\in F^2_{pot}$ and $N>0$ ($N\in \N^*$ in the discrete case)
\begin{equation}\label{dssskdjdjdhjhd090334}
\|\Pi_N \big(\xi_{pot}(\eta)\big)- \big(\Pi_N \xi_{pot}(\eta)\big)_{pot}\|_{\L^2(T^d_N)}\leq C_d \big<\xi_{pot}^2\big>^\frac{1}{2}\big(f(N,\xi_{pot},\eta)\big)^{-\frac{1}{2}}
\end{equation}
The same rate of convergence can be obtained with the second strategy in the discrete case for both equations \eref{dkdjdjdhjhd090334} and \eref{dkdjdjdhjhd09033445}. Then following the proof of our applications one can relate the rate at which the effective conductivity can be approximated by its periodizations to the function
$f(N,\xi_{pot},\eta)$.
\end{Remark}

\subsubsection{Proof of theorem \ref{jhdchdcbcc221}}
We will now give the proof of theorem \ref{jhdchdcbcc221} based on proposition \ref{dddhdhdhhh8711}. We will first give the proof in the continuous case. First, let us remind the standard ergodic theorem that we will use.\\
Let $f(x)\in L^1_{loc}(\R^d)$. A number $M\{f\}$ is called the mean value of $f$ if 
\begin{equation}\label{sajhvs}
\lim_{\epsilon \rightarrow 0} \int_K f(\epsilon^{-1}x)dx=|K| M\{f\}
\end{equation}
for any Lebesgue measurable bounded set $K\subset \R^d$ (here $|K|$ stands for the Lebesgue measure of $K$). Let $K_t=\{x \in \R^d, \, t^{-1}x \in K\}$ denote the homothetic dilatation, with ratio $t>0$, of the set $K$. Then
\eref{sajhvs} can be written in a more habitual form:
\begin{equation}\label{sajhvs2}
\lim_{t\rightarrow \infty}\frac{1}{t^d |K|}\int_{K_t}f(x)dx=M\{f\}
\end{equation}
The following theorem is the theorem 10 of the chapter VIII.7.10 of \cite{DuScha67} (see also the theorem 7.2 of \cite{JiKoOl91} ).
\begin{Theorem}\label{wlkkiudbu1}
Let $f\in L^p(\mu)$ with $1\leq p<\infty$. Then for almost all $\eta \in X$ the realization $T_x f(\eta)$ posses a mean value in the sense of \eref{sajhvs2}. Moreover, the mean value $M\{T_xf(\eta)\}$, considered as a function of $\eta \in X$ is invariant, and for almost all $\eta \in X$
\begin{equation}
<f>\equiv \int_X f(\eta)\mu(d\eta)=M\{T_x f(\eta)\}
\end{equation}
If $p>1$, then the limit in \eref{sajhvs2} also exists in the norm of $L^p$ and the functions are for $t>0$, all dominated by a function in $L^p$.
\end{Theorem}
Now let us observe that from equation \ref{dkdjdjdhjhd090334} one easily obtains that 
\begin{equation}\label{aadjhfgfgf8817222}
F^2_{pot}\subset \lim_{N\rightarrow \infty }F^2_{pot}(T^d_N)
\end{equation}
similarly from equation \ref{dkdjdjdhjhd09033445} one easily obtains that
\begin{equation}\label{aaaddddddhhjfhfhfh8811}
F^2_{sol}\subset \lim_{N\rightarrow \infty }F^2_{sol}(T^d_N)
\end{equation}
Now let  $\xi \in \lim_{N\rightarrow \infty }F^2_{pot}(T^d_N)$ and $\nu \in \L^2_{sol}(T^d_N)$.
Using the ergodic theorem \ref{wlkkiudbu1}, one has $\mu$-a.s.
\begin{equation}
\lim_{N\rightarrow \infty}N^{-d}\int_{T^d_N}\Pi_N \xi(x,\eta).\Pi_N \nu(x,\eta)\,dx=\big<\xi.\nu\big>
\end{equation}
 Using
\begin{equation}
\begin{split}
\int_{T^d_N}\Pi_N \xi(x,\eta).\Pi_N \nu(x,\eta)\,dx&=\int_{T^d_N}\big(\Pi_N \xi(x,\eta)\big)_{pot}.\big(\Pi_N \nu(x,\eta)\big)_{pot}\,dx\\&+
\int_{T^d_N}\Big(\Pi_N \xi(\eta)-\big(\Pi_N \xi(\eta)\big)_{pot}\Big).\Pi_N \nu(x,\eta)\,dx
\end{split}
\end{equation}
it follows that
\begin{equation}
\begin{split}
\big|N^{-d}\int_{T^d_N}\Pi_N \xi(x,\eta).&\Pi_N \nu(x,\eta)\,dx\big|\leq \|\big(\Pi_N \xi(\eta)\big)_{pot}\|_{\L^2(T^d_N)}
\|\big(\Pi_N \nu(\eta)\big)_{pot}\|_{\L^2(T^d_N)}\\&+
\|\Pi_N \big(\xi(\eta)\big)- \big(\Pi_N \xi(\eta)\big)_{pot}\|_{\L^2(T^d_N)}\|\big(\Pi_N \nu(\eta)\big)\|_{\L^2(T^d_N)}
\end{split}
\end{equation}
Using the ergodic theorem \ref{wlkkiudbu1} and proposition \ref{dddhdhdhhh8711}, one easily obtains that $\mu$-a.s.
\begin{equation}
\lim_{N\rightarrow \infty}N^{-d}\int_{T^d_N}\Pi_N \xi(x,\eta).\Pi_N \nu(x,\eta)\,dx=0
\end{equation}
which proves that $\big<\xi.\nu\big>=0$. Thus $\xi \bot \L^2_{sol}(X,\mu)$, which implies from Weyl decomposition that $\xi \in F^2_{pot}(X,\mu)$. Thus we have proven that 
\begin{equation}\label{aadjhfgfgfsss88172s22}
\lim_{N\rightarrow \infty }F^2_{pot}(T^d_N)\subset F^2_{pot} 
\end{equation}
Similarly one proves that
\begin{equation}\label{aadjhfgfgfsss881ss7222}
\lim_{N\rightarrow \infty }F^2_{sol}(T^d_N)\subset F^2_{sol} 
\end{equation}
Combining the equations \eref{aadjhfgfgf8817222}, \eref{aaaddddddhhjfhfhfh8811}, \eref{aadjhfgfgfsss88172s22} and \eref{aadjhfgfgfsss881ss7222} one concludes the poof of theorem \ref{jhdchdcbcc221}.
The proof in the discrete case being similar, we will just remind below the ergodic theorem which is used.
For any  bounded set $K\subset \Z^d$. Let $K_t=\{x \in \Z^d, \, t^{-1}x \in K\}$ denote the homothetic dilatation, with ratio $t>0$, of the set $K$.
Let $f(x)\in L^1_{loc}(\Z^d)$. A number $M\{f\}$ is called the mean value of $f$ if 
\begin{equation}\label{sajhvs002}
\lim_{t\rightarrow \infty}\frac{1}{t^d |K|}\sum_{x\in K_t}f(x)=M\{f\}
\end{equation}
For any  bounded set $K\subset \Z^d$.\\ 
The following theorem is the theorem 9 of the chapter VIII.6.9 of \cite{DuScha67} 
\begin{Theorem}\label{wlkkiudbu001}
Let $f\in L^p(\mu)$ with $1\leq p<\infty$. Then for almost all $\eta \in X$ the realization $T_x f(\eta)$ posses a mean value in the sense of \eref{sajhvs002}. Moreover, the mean value $M\{T_xf(\eta)\}$, considered as a function of $\eta \in X$ is invariant, and for almost all $\eta \in X$
\begin{equation}
<f>\equiv \int_X f(\eta)\mu(d\eta)=M\{T_x f(\eta)\}
\end{equation}
If $p>1$, the limit in \eref{sajhvs002} also exists in the norm of $L^p$ and the functions are for $t>0$, all dominated by a function in $L^p$.
\end{Theorem}

\subsection{Periodic approximation of variational functionals on the ergodic spaces}
In this subsection it will be shown that the almost sure continuity of Weyl's decomposition has a direct consequence on variational functionals on the ergodic space. The results and notations given in this subsection are valid in the continuous case as in the discrete case.
Let $m,p\in \N^2$, and 
\begin{equation}
\Phi:\quad (\eta,X^1,\ldots,X^m,Y^1,\ldots,Y^p)\longrightarrow \Phi(\eta,X^1,\ldots,X^m,Y^1,\ldots,Y^p)
\end{equation}
a mapping from $X\times \R^{d\times m}\times \R^{d\times p}$ into $\R^+$
such that for any \\$(\xi^1,\ldots,\xi^m,\nu^1,\ldots,\nu^p)\in (F^2_{pot})^m\times (F^2_{sol})^p$ one has
$$\eta \rightarrow\Phi(\eta,\xi^1,\ldots,\xi^m,\nu^1,\ldots,\nu^p)\in \L^1(X,\mu).$$
We write 
\begin{equation}\label{dlkjnjdnjdnnjnvv8831}
Z(\Phi,\mu):=\inf_{(\xi^1,\ldots,\xi^m,\nu^1,\ldots,\nu^p)\in (F^2_{pot})^m\times (F^2_{sol})^p}\Big<\Phi(\eta,\xi^1,\ldots,\xi^m,\nu^1,\ldots,\nu^p)\Big>
\end{equation}
Let us define  for $N>0$ the function (in the continuous case)
\begin{equation}
\begin{split}
\Psi_N: X\times(\L^2(T^d_N))^m\times (\L^2(T^d_N))^p \rightarrow& \R^+\\
\big(\eta,v^1,\ldots,v^m,q^1,\ldots,q^p\big)\rightarrow&  N^{-d}
\int_{x\in [0,N(^d}\Phi(\tau_{-x}\eta,v^1(x),\ldots,v^m(x)\\&\qquad,q^1(x),\ldots,q^p(x))\,dx
\end{split}
\end{equation}
In the discrete case we shall write
\begin{equation}
\begin{split}
\Psi_N: X\times(\L^2(T^d_N))^m\times (\L^2(T^d_N))^p \rightarrow& \R^+\\
\big(\eta,v^1,\ldots,v^m,q^1,\ldots,q^p\big)\rightarrow&  N^{-d}
\sum_{x\in \Z\cap [0,N(^d}\Phi(\tau_{-x}\eta,v^1(x),\ldots,v^m(x)\\&\qquad,q^1(x),\ldots,q^p(x))
\end{split}
\end{equation}
Let us define for $N\in \N$ the random variable $Z(N,\eta)$ by
\begin{equation}\label{jshdgddiudgg8981221}
\begin{split}
Z(N,\eta):=\inf_{(v^1,\ldots,v^m,q^1,\ldots,q^p)\in (F^2_{pot}(T^d_N))^m\times (F^2_{sol}(T^d_N))^p}
\Psi_N\big(\eta,v^1,\ldots,v^m,q^1,\ldots,q^p\big)
\end{split}
\end{equation}
Observe that $Z(N,\eta)$ corresponds to the periodization of the variational problem associated to $Z(\Phi,\mu)$ over the torus $T^d_N$ for a particular realization $\eta$ of the ergodic space.
\begin{Definition}
We say that the function $\Phi$ is admissible if there exists a strictly increasing continuous function $g$ from $\R^+$ into $\R^+$ such that for all $\eta\in X$, the function 
\begin{equation}
\begin{split}
g\circ\Psi_N(\eta): (\L^2(T^d_N))^m\times (\L^2(T^d_N))^p \rightarrow& \R^+\\
\big(v^1,\ldots,v^m,q^1,\ldots,q^p\big)\rightarrow&  g\Big(\Psi_N\big(\eta,v^1,\ldots,v^m,q^1,\ldots,q^p\big)\Big)
\end{split}
\end{equation}
 is upper semi-continuous with respect to the norm $\sum_{i=1}^m\|v^i\|_{\L^2(T^d_N)}+\sum_{j=1}^p\|v^j\|_{\L^2(T^d_N)}$ a.s. uniformly in $N$ and $\eta$.
\end{Definition}
\begin{Theorem}\label{sddhdvdhdzhvdzd78712}
If the function $\Phi$ is admissible then for $\mu$-almost all $\eta\in X$
\begin{equation}
\lim \sup_{N\rightarrow \infty} Z(N,\eta)\leq Z(\Phi,\mu)
\end{equation}
\end{Theorem}
\begin{proof}
We will write the theorem for the continuous case. In the discrete case the proof is trivially similar.
Let $(\xi^1,\ldots,\xi^m,\nu^1,\ldots,\nu^p)\in (F^2_{pot})^m\times (F^2_{sol})^p$. To prove the theorem it is sufficient to show that for $\mu$-almost all $\eta\in X$ 
\begin{equation}
\lim \sup_{N\rightarrow \infty} g\big(Z(N,\eta)\big)\leq g\Big(\Big<\Phi(\eta,\xi^1,\ldots,\xi^m,\nu^1,\ldots,\nu^p)\Big>\Big)
\end{equation}
By the equation \eref{jshdgddiudgg8981221} we have
\begin{equation}\label{jsaaahdgddiudgg89ccwq81221}
\begin{split}
g\big(Z(N,\eta)\big)= J_1(N)+J_2(N)
\end{split}
\end{equation}
with 
\begin{equation}\label{jshdgddiudgg89ccwq81221aa}
\begin{split}
J_1(N)= g\Big(N^{-d}\int_{x\in [0,N(^d}\Phi(\tau_{-x}\eta,\xi^1(x,\eta),\ldots,\xi^m(x,\eta),\nu^1(x,\eta),\ldots,\nu^p(x,\eta))\,dx\Big)
\end{split}
\end{equation}
and
\begin{equation}\label{jshdgddiudgg89ccwq81221bb}
\begin{split}
J_2(N)=& \inf_{(v^1,\ldots,v^m,q^1,\ldots,q^p)\in (F^2_{pot}(T^d_N))^m\times (F^2_{sol}(T^d_N))^p}
\Big(g\Big(\Psi_N\big(\eta,v^1,\ldots,v^m,q^1,\ldots,q^p\big)\Big)\\&-g\Big(\Psi_N\big(\eta,\Pi_N\xi^1,\ldots,\Pi^N\xi^m,\Pi^N \nu^1,\ldots,\Pi^N \nu^p \big)\Big)\Big)
\end{split}
\end{equation}
by the ergodic theorem \ref{wlkkiudbu1}  one has for $\mu$-almost all $\eta\in X$
\begin{equation}\label{jshdgddiudgg89ccwq81221assa}
\begin{split}
\lim_{N\rightarrow \infty}J_1(N)= g\Big(\Big<\Phi(\eta,\xi^1,\ldots,\xi^m,\nu^1,\ldots,\nu^p)\Big>\big)
\end{split}
\end{equation}
Now by  theorem  \ref{jhdchdcbcc221} for $\mu$-almost all $\eta\in X$ there exists  sequences\\ $v^{1,N},\ldots, q^{1,N},v^{1,N},\ldots, q^{p,N}\in \big(F^2_{pot}(T^d_N)\big)^m\times \big(F^2_{sol}(T^d_N)\big)^p$ such that
\begin{equation}\label{gddhgdgfdfff77611}
\lim_{N\rightarrow \infty} \sum_{i=1}^m \|\Pi_N \xi^i(\eta)- v^{i,N}\|_{\L^2(T^d_N)}
+\sum_{j=1}^p \|\Pi_N \nu^j(\eta)- v^{j,N}\|_{\L^2(T^d_N)}=0
\end{equation}
From equation \eref{jshdgddiudgg89ccwq81221bb} one obtains
\begin{equation}\label{jshdgddiussadgg89ccwq81221bb}
\begin{split}
J_2(N)\leq &g\Big(\Psi_N\big(\eta,v^{1,N},\ldots,v^{m,N},q^{1,N},\ldots,q^{p,N}\big)\Big)\\&-g\Big(\Psi_N\big(\eta,\Pi_N\xi^1,\ldots,\Pi^N\xi^m,\Pi^N \nu^1,\ldots,\Pi^N \nu^p \big)\Big)
\end{split}
\end{equation}
Combining \eref{jshdgddiussadgg89ccwq81221bb} with \eref{gddhgdgfdfff77611} and the uniform upper semi-continuity of $g\circ \psi(\eta)$ we concludes that for $\mu$-almost all $\eta\in X$
\begin{equation}\label{jshdgddiussadgg89ccwq81sss221bb}
\begin{split}
\lim \sup_{N\rightarrow \infty} J_2(N)\leq 0
\end{split}
\end{equation}
Which concludes the proof of the theorem.
\end{proof}

\section{Application}
\subsection{Symmetric continuous Operator}\label{KKKKKKKK00001}
\subsubsection{Homogenization in the ergodic medium}
Let $A(\eta)$ be a $d\times d$ bounded symmetric matrix defined on $X$ ($A_{i,j}\in L^\infty(X,\mu)$) and satisfying the following ellipticity condition
\begin{equation}\label{ddjhbhdh8913}
\nu_1 |\xi|^2 \leq ^t\xi A \xi\leq \nu_2 |\xi|^2,\quad \nu_1 >0
\end{equation}
for almost all $\eta\in X$. Realizations $A(x,\eta)=T_x A(\eta)$ of this matrix are considered and we are interested in describing the homogenization for almost all $\eta \in X$ of the operator $\nabla_x A(x,\eta)\nabla_x$.\\
Consider $\sigma(A,\mu)$ the $d\times d$ positive definite symmetric matrix defined by the following variational formula: for $\xi\in\R^d$
\begin{equation}\label{dkjhdblkclkbu2}
^t\xi \sigma(A,\mu)\xi=\inf_{v\in F^2_{pot}}\big<{^t(}\xi+v)A(\xi+v)\big>
\end{equation}
Observe that $\sigma(A,\mu)$ corresponds to the effective conductivity associated to the operator $\nabla_x A(x,\eta)\nabla_x$. Indeed by the theorem 7.4 of \cite{JiKoOl91}
 for any bounded domain $Q\subset \R^d$ and any $f\in H^{-1}(Q)$ the solutions $u^\epsilon$ of the Dirichlet problems ($A^\epsilon(x,\eta)=A(x/\epsilon,\eta)$)
\begin{equation}\label{dhsdkhkbc0}
\nabla A^\epsilon \nabla u^\epsilon=f,\quad u^\epsilon \in H^1_0(Q)
\end{equation}
possess the following properties of convergence (in the weak topology)
\begin{equation}\label{dhsdkhkbc1}
u^\epsilon\rightarrow u^0 \;\text{in}\; H^1_0(Q),\quad A^\epsilon \nabla u^\epsilon \rightarrow \sigma(A,\mu) \nabla u^0
\;\text{in}\; \L^2(Q)
\end{equation}
where $u^0$ is the solution of the Dirichlet problem
\begin{equation}\label{dhsdkhkbc2}
\nabla \sigma(A,\mu) \nabla u^0=f,\quad u^0 \in H^1_0(Q)
\end{equation}
Moreover writing $y_t^\eta$ the diffusion associated to the operator $\nabla_x A(x,\eta)\nabla_x$, 
and $\P_\eta$ the law of that started from $0$ in $\R^d$ it is well known (\cite{KiVa86}, \cite{JiKoOl91}, \cite{Ol94}) that under the law $\mu \otimes \P_\eta$,
$\epsilon y_{t/\epsilon^2}^\eta$ converges in law as $\epsilon \downarrow 0$ towards a Brownian Motion starting from $0$ with covariance matrix (effective diffusivity) $2 \sigma (A)$.

\subsubsection{Periodization of the ergodic medium}
For $\eta \in X$, we write $A^N(x,\eta)$ obtained by periodizing $A(x,\eta)$ over the torus $T^d_N$ (of dimension $d$ and side $N$, $\R^d\big/ (N\Z^d)$)
\begin{equation}
A^N(x,\eta)=A(x-N[x/N],\eta)
\end{equation}
where $[y]$ is the integer part of $y$. For $\eta \in X$, we define $\sigma(A^N,\eta)$ the $d\times d$ symmetric positive definite matrix by the  following variational formula: for $\xi\in\R^d$
\begin{equation}\label{ddvjvdjhvh881}
^t\xi \sigma(A^N,\eta)\xi=\inf_{f\in C^\infty(T^d_N)}N^{-d}\int_{T^d_N}{^t(}\xi+\nabla f(x))A^N(x,\eta)(\xi+\nabla f(x))\,dx
\end{equation}
Observe that $\sigma(A^N,\eta)$ corresponds to the effective conductivity associated to the periodic operator $\nabla A^N(x,\eta)\nabla$ in the sense given above in the equations \eref{dhsdkhkbc0}, \eref{dhsdkhkbc1} and \eref{dhsdkhkbc1}. Writing $y_t^{\eta,N}$ the diffusion associated to the operator $\nabla_x A^N(x,\eta)\nabla_x$, it is well known (\cite{JiKoOl91}, \cite{Ol94}) that 
$\epsilon y_{t/\epsilon^2}^{\eta,N}$ converges in law as $\epsilon \downarrow 0$ towards a Brownian Motion starting from $0$ with covariance matrix (effective diffusivity) $2 \sigma \big(A^N,\eta\big)$. Notice, whereas $\sigma(A,\mu)$ is a constant (not random) matrix,  $\sigma \big(A^N,\eta\big)$ is a random matrix on $X$, which depends on the particular realization $A^N(x,\eta)$ of the periodic environment.

\subsubsection{The main theorem}
It is our purpose to prove the following theorem
\begin{Theorem}\label{JHlkhbsljbj81}
For $\mu$-almost all $\eta \in X$
\begin{equation}
\lim_{N\rightarrow +\infty}\sigma(A^N,\eta)=\sigma(A,\mu)
\end{equation}
\end{Theorem}

\subsubsection{Proof}
Let $\xi \in \R^d$.
Let us apply theorem \ref{JHlkhbsljbj81} with $m=1$, $p=0$ and $\Phi(\eta, X^1)={^t(\xi+X^1)A(\eta)(\xi+X^1)}$. By  the  Minkowski inequality and the uniform ellipticity condition \eref{ddjhbhdh8913} one has that for $N\in \N^*$ and $\mu$-almost all $\eta \in X$, $v^1,v^2 \in \big(\L^2(T^d_N)\big)^2$
\begin{equation}
\big(\Psi_N(\eta,v^1)\big)^{\frac{1}{2}}-\big(\Psi_N(\eta,v^2)\big)^{\frac{1}{2}}\leq \nu_2 \|v^1-v^2\|_{\L^2(T^d_N)}
\end{equation}
It follows that $\Phi$ is admissible and from the variational formulae
 \eref{ddvjvdjhvh881}, \eref{ddvjvdjhvh881} and theorem \ref{JHlkhbsljbj81} one obtains that for $\mu$-almost all $\eta \in X$
\begin{equation}
\lim\sup_{N\rightarrow \infty}{^t\xi} \sigma(A^N,\eta)\xi\leq {^t\xi}\sigma(A,\mu)\xi
\end{equation}
Which gives the upper bound of theorem \ref{JHlkhbsljbj81}. For the lower bound we will  apply theorem \ref{JHlkhbsljbj81} with $m=0$, $p=1$ and $\Phi(\eta, Y^1)={^t(\xi+Y^1)A^{-1}(\eta)(\xi+Y^1)}$.
By the   Minkowski inequality and the uniform ellipticity condition \eref{ddjhbhdh8913} one has that for $N\in \N^*$ and $\mu$-almost all $\eta \in X$, $q^1,q^2 \in \big(\L^2(T^d_N)\big)^2$
\begin{equation}
\big(\Psi_N(\eta,q^1)\big)^{\frac{1}{2}}-\big(\Psi_N(\eta,q^2)\big)^{\frac{1}{2}}\leq (\nu_1)^{-1} \|q^1-q^2\|_{\L^2(T^d_N)}
\end{equation}
Moreover let us remind the
following well known (\cite{JiKoOl91}) variational formulas: for $l\in \R^d$
\begin{equation}\label{shsjshbsjhb801}
^tl \sigma(A,\mu)^{-1}l=\inf_{p\in F^2_{sol}}\big<(l+p)A^{-1}(l+p)\big>
\end{equation}
\begin{equation}\label{shsjshbsjhb81}
^tl \sigma(A^N,\eta)^{-1}l=\inf_{\nu \in F^2_{sol}(T^d_N)}N^{-d}\int_{T^d_N}{^t\big(}l+\nu(x)\big)\big(A^N(x,\eta)\big)^{-1}\big(l+\nu(x)\big)\,dx
\end{equation}
Then it follows that $\Phi$ is admissible and from the variational formulae
 \eref{shsjshbsjhb801}, \eref{shsjshbsjhb81} and theorem \ref{JHlkhbsljbj81} one obtains that for $\mu$-almost all $\eta \in X$
\begin{equation}
\lim\sup_{N\rightarrow \infty}{^tl} \big(\sigma(A^N,\eta)\big)^{-1}l\leq {^tl} \big(\sigma(A,\mu)\big)^{-1}l
\end{equation}
Which gives the lower  bound of theorem \ref{JHlkhbsljbj81}.

\subsection{Non symmetric continuous Operator, diffusion in divergence free flow}\label{KKKKKKKK00002}
\subsubsection{Homogenization in the ergodic medium}
Let $E$ be a  $d\times d$ bounded skew-symmetric matrix defined on $X$ ($E_{i,j}\in L^\infty(X,\mu)$). Let $a$ be a constant symmetric positive definite $d\times d$ matrix.
Realizations $E(x,\eta)=T_x E(\eta)$ of this matrix are considered and we are interested in describing the homogenization for almost all $\eta \in X$ of the operator 
\begin{equation}
L_{E}=\nabla_x \big(a+E(x,\eta)\big)\nabla_x
\end{equation}
$E$ is seen as the stream matrix of the incompressible flow $^t\nabla.E$.\\
Let $z_t^\eta$ be the process generated by $L_E$,  
and $\P_\eta$ the law of that diffusion started from $0$ in $\R^d$. It is well known (see for instance \cite{Ol94})that under the law $\mu \otimes \P_\eta$ as $\epsilon \downarrow 0$, $\epsilon z_{t/{\epsilon^2}}^\eta$ converges in law to a Brownian motion with covariance matrix $D(a,E,\mu)$: for $l \in \R^d$
\begin{equation}
^tlD(a,E,\mu)l=2{^tl}al+2\big<|v_l|_a^2\big>
\end{equation}
Where we have used the notation $|\xi|^2_a:={^t\xi}a\xi$ for $\xi\in \R^d$ and $v_l$ defined as the unique solution $u\in F^2_{pot}$ of 
\begin{equation}\label{wlbuasssbu1}
<\phi.(a+E) (l+u)>=0,\quad \forall \phi \in F^2_{pot};\quad u\in F^2_{pot}
\end{equation}
The existence of a solution for this problem follows from the Lax-Milgram Lemma and the estimate $<v.(a+E)v>\geq \lambda_{\min}(a) \|v\|^2_{\L^2(X,\mu)}$; \cite{JiKoOl91}.\\
Obviously, the solution $v_l$ of the problem \eref{wlbuasssbu1} depends linearly on $l \in \R^d$. Therefore \\$\big<(a+E)(l+v_l)\big>$ is a linear form with respect to $l$. The effective conductivity $\sigma(a,E,\mu)$ is defined by
\begin{equation}\label{dshdjh88171}
\sigma(a,E,\mu)l=\big<(a+E)(l+v_l)\big>
\end{equation}
It  is a non-symmetric matrix relating the gradient of the heat intensity with the flux \cite{FaPa94} by \eref{dshdjh88171}.
Observe that the symmetric part of the  effective conductivity gives the effective diffusivity by the following relation:
\begin{equation}
D(a,E,\mu)=2 \sigma_\sym(a,E)
\end{equation}

\subsubsection{Periodization of the ergodic medium}
For $\eta \in X$, we write $E^N(x,\eta)$ obtained by periodizing $E(x,\eta)$ over the torus $T^d_N$ 
\begin{equation}
E^N(x,\eta)=E(x-N[x/N],\eta)
\end{equation}
We are interested in describing the homogenization for almost all $\eta \in X$ of the operator 
\begin{equation}
L_{E}^N=\nabla_x \big(a+E^N(x,\eta)\big)\nabla_x
\end{equation}
Let $z_t^{\eta,N}$ be the process generated by $L_E^N$. It is well known (see for instance \cite{Nor97})that as $\epsilon \downarrow 0$, $\epsilon z_{t/{\epsilon^2}}^{\eta,N}$ converges in law to a Brownian motion with covariance matrix $D(a,E^N,\eta)$ with for $l \in \R^d$
\begin{equation}
^tlD(a,E^N,\eta)l=2{^tl}al+2N^{-d}\int_{T^d_N}|\psi_l(x,\eta)|_a^2\,dx
\end{equation}
Where  $\psi_l$ defined as the unique solution $\psi \in H^1(T^d_N)$  of 
\begin{equation}\label{wlbuasssbusqw21}
\int_{T^d_N} \phi(x)\big(a+E^N(x,\eta)\big) \big(l+\psi(x)\big)\,dx=0,\quad \forall \phi \in H^1(T^d_N);\quad \psi\in H^1(T^d_N)
\end{equation}
We have noted $H^1(T^d_N)$ the closure of $\{\nabla f\,:\, f\in C^\infty(T^d_N)\}$ in $L^2(T^d_N)$ with respect to the $L^2$-norm.
Obviously, the solution $\psi_l$ of the problem \eref{wlbuasssbusqw21} depends linearly on $l \in \R^d$. Therefore $\int_{T^d_N}\big(a+E^{\eta,N}(x)\big)\big(l+\psi_l(x,\eta)\big)\,dx$ is a linear form with respect to $l$. The effective conductivity $\sigma(a,E^N,\eta)$ is defined by: for $l\in \R^d$
\begin{equation}\label{dshdjh88170001}
\sigma(a,E^N,\eta)l=\int_{T^d_N}\big(a+E^{\eta,N}(x)\big)\big(l+\psi_l(x)\big)\,dx
\end{equation}
It  is a non-symmetric matrix relating the gradient of the heat intensity with the flux \cite{FaPa94} by \eref{dshdjh88170001}.
Observe that the symmetric part of the  effective conductivity gives the effective diffusivity by the following relation:
\begin{equation}
D(a,E^N,\eta)=2 \sigma_{\sym}(a,E^N,\eta)
\end{equation}

\subsubsection{The main theorem}
It is our purpose to prove the following theorem
\begin{Theorem}\label{dkhbdlhl8981}
For $\mu$-almost all $\eta \in X$
\begin{equation}
\lim_{N\rightarrow +\infty}\sigma(a,E^N,\eta)=\sigma(a,E,\mu)
\end{equation}
In particular
\begin{equation}
\lim_{N\rightarrow +\infty}D(a,E^N,\eta)=D(a,E,\mu)
\end{equation}
\end{Theorem}

\subsubsection{Core of the proof: Variational formulations and theorem \ref{JHlkhbsljbj81}}
As for a symmetric operator,  the proof theorem \ref{dkhbdlhl8981} relies theorem \ref{jhdchdcbcc221} and the variational formulae associated to the effective conductivity.
\paragraph{Variational Formulation of the effective diffusivity in the periodic case}
For $N>0$, let us write $\S(T^d_N)$  the set of  skew-symmetric matrices with smooth coefficients defined on $T^d_N$ and for $H\in \S(T^d_N)$, $\Div H$ is the vector field defined by $(\Div H)_i= \sum_{j=1}^d \partial_j H_{i,j}$.\\
In the periodic case, we will use Norris's variational formulation (obtained by polarization \cite{Nor97}) to  control $\sigma(a,E^N,\eta)$. For $y\in \R^d$ we will write $|y|^2_{a^{-1}}:={^ty a^{-1}y}$.\\
For all $\xi,l \in \R^d$,
\begin{equation}\label{dsjkdkjsdh8171}
\begin{split}
|\xi-\sigma(a,E^N,\eta) l&|_{\sigma_{\sym}^{-1}(a,E^N,\eta)}^2=\\&\inf_{f, H\in C^\infty(T^d_N)\times \S(T^d_N)} N^{-d}\int_{T^d_N}|\xi-\nabla H-(a+E^N(x,\eta))(l-\nabla f)|_{a^{-1}}^2\,dx
\end{split}
\end{equation}
For all $l\in \R^d$
\begin{equation}\label{HoCoeDNorPre14}
|l|_{\sigma_{\sym}(a,E^N,\eta)}^2=\inf_{\xi \perp l, f,H \in C^\infty(T^d_N)\times \S(T^d_N)}N^{-d}\int_{T^d_N}|\xi-\nabla H-(a+E^N(x,\eta))(l-\nabla f)|^2_{a^{-1}}\,dx
\end{equation}
Where we have written $\xi \perp l:=\{\xi\in \R^d\,:\,\xi.l=0\}$. We also have for all $\xi\in\R^d$
\begin{equation}\label{HoCoeDNorPre12}
|\xi|_{\sigma_{\sym}^{-1}(a,E^N,\eta)}^2=\inf_{f,H\in C^\infty(T^d_N)\times \S(T^d_N)}N^{-d}\int_{T^d_N}|\xi-\nabla H+(a+E^N(x,\eta))\nabla f|_{a^{-1}}^2\,dx
\end{equation}
\begin{Remark}
Let us remind, as it has been noticed by J.R. Norris \cite{Nor97}, that from \eref{HoCoeDNorPre14} and \eref{HoCoeDNorPre12} one obtains that
\begin{equation}\label{ksjdkjdbj1}
a \leq \sigma_{\sym}(a,E^N,\eta)\leq a+ N^{-d}\int_{T^d_N}{^tE^N(x,\eta)}a^{-1}E^N(x,\eta)\,dx
\end{equation}
For a saddle point variational formulation we refer to \cite{FaPa94}.
\end{Remark}

\paragraph{ Variational Formulation of the effective diffusivity in the ergodic case}
The following theorem proven in subsection \ref{hdodkdbbii871} is  inspired from the variational formulation given for the periodic case by J. R. Norris \cite{Nor97} (lemma 3.1), (for a non local variational formulation we refer to \cite{FaPa96})
\begin{Theorem}\label{dkhdbldbdhjb81}
For all $\xi,l \in \R^d$,
\begin{equation}\label{dsdjhdldhi8101}
\begin{split}
|\xi-\sigma(a,E,\mu) l|_{\sigma_{\sym}^{-1}(a,E,\mu)}^2=\inf_{v,p\in F_{pot}\times F_{sol}}\Big<|\xi-p-(a+E)(l-v)|_{a^{-1}}^2\Big>
\end{split}
\end{equation}
for $l \in\R^d$
\begin{equation}\label{ddlkjdj1}
^tl \sigma_{\sym}(a,E,\mu)l=\inf_{\xi \perp l,v\in F^2_{pot},p\in F^2_{sol}}\big<|\xi-p-(a+E)(l-v)|^2_{a^{-1}}\big>
\end{equation}
For all $\xi \in \R^d$
\begin{equation}\label{HoCoeDNorPdddre12}
|\xi|_{\sigma_{\sym}^{-1}(a,E,\mu)}^2=\inf_{v,p\in F_{pot}\times F_{sol}}\Big<|\xi-p+(a+E)v|_{a^{-1}}^2\Big>
\end{equation}
\end{Theorem}
\begin{Remark}
Let us observe that from \eref{ddlkjdj1} and \eref{HoCoeDNorPdddre12} one obtains that
\begin{equation}\label{ksjdkjdbj2}
a \leq \sigma_{\sym}(a,E,\mu)\leq a+ \big<^tEa^{-1}E\big>\end{equation}
\end{Remark}
Let $\xi,l \in \R^d$.
Let us apply theorem \ref{JHlkhbsljbj81} with $m=1$, $p=1$ and $\Phi(\eta, X^1,Y^1)=|\xi-Y^1-(a+E)(l-X^1)|_{a^{-1}}^2$. By  the  Minkowski inequality one has that for $N\in \N^*$ and $\mu$-almost all $\eta \in X$, $v^1,v^2,q^1,q^2 \in \big(\L^2(T^d_N)\big)^4$
\begin{equation}
\begin{split}
\big(\Psi_N(\eta,v^1,q^1)\big)^{\frac{1}{2}}&-\big(\Psi_N(\eta,v^2,q^2)\big)^{\frac{1}{2}}\leq
C_d \big(\lambda_{\min}(a)\big)^{-1/2} \|q^1-q^2\|_{\L^2(T^d_N)}\\&+C_d (\lambda_{\max}(a)+\|E\|_{L^\infty(X,\mu)})\big(\lambda_{\min}(a)\big)^{-1/2}\|v^1-v^2\|_{\L^2(T^d_N)}
\end{split}
\end{equation}
It follows that $\Phi$ is admissible and from the variational formulae
 \eref{dsdjhdldhi8101}, \eref{dsjkdkjsdh8171} and theorem \ref{JHlkhbsljbj81} one obtains that for $\mu$-almost all $\eta \in X$
\begin{equation}
\lim\sup_{N\rightarrow \infty}|\xi-\sigma(a,E^N(\eta) ) l|_{\sigma_{\sym}^{-1}(a,E^N,\eta)}^2\leq |\xi-\sigma(a,E,\mu) l|_{\sigma_{\sym}^{-1}(a,E,\mu)}^2
\end{equation}
Choosing $\xi:=\sigma(a,E,\mu) l$ in this equation, one obtains from \eref{ksjdkjdbj1} and \eref{ksjdkjdbj2} that for $\mu$-almost all $\eta\in X$
\begin{equation}
\lim_{N\rightarrow \infty}\sigma(a,E^N(\eta) ) l =\sigma(a,E,\mu) l
\end{equation} 
Which concludes the proof of theorem \ref{dkhbdlhl8981}.

\subsection{Discrete Operator}\label{KKKKKKKK00003}
We shall extend in this subsection our results to the discrete case.  
\subsubsection{The ergodic homogenization problem}
We will consider a symmetric random walk on $\Z^d$ as in \cite{CapIof01} but with ergodic jump rates instead of i.i.d. 
The random ergodic environment will be represented by the random $d$-dimensional vector $\xi_i(\eta)$ ($i\in \{1,\ldots,d\}$) on $X$, we will write $\xi_i(x,\eta)=\xi_i(\tau_{-x}\eta)$. We will assume that there exists $c\geq 1$ such that for $\mu$-almost all $\eta\in X$, 
\begin{equation}\label{djhdjhwhdheu1}
1/c\leq \xi_i(\eta)\leq c
\end{equation}
Let us write $X\big(t,\xi(\eta)\big)$ the nearest neighbor symmetric random walk on $\Z^d$ with jump according to $\xi_i(x,\eta)$ rates ($\xi_i(x,\eta)$ is the jump rate between from the site $x$ to the site $x+e_i$ and also from the site $x+e_i$ to the site $x$).\\ In the quenched regime (for a fixed $\eta$), $\P^{\xi(\eta)}_x$ stands for the probability law of this process when the walk starts at $x\in \Z^d$. It is well known (\cite{KiVa86}, \cite{MaFe89}, \cite{CapIof01}) that in the annealed regime (under the law $\mu\otimes \P^{\xi(\eta)}_0$) as $\epsilon \downarrow 0$, $\epsilon X\big(t/\epsilon^2,\xi(\eta)\big)$ converges in law towards a Brownian Motion  with covariance matrix (effective diffusivity) $D(\xi,\mu)$.

\subsubsection{Periodization of the ergodic medium}
For $N\in \N^*$ and $\eta \in X$ we write $\xi^N(\eta)$ the periodized bond configuration associated to $\xi(x,\eta)$ over the torus $T_N=\Z^d/ N \Z^d$. For $x\in \Z^d$ decomposed as $x=y+N z$ with $y \in \{0,\ldots,N-1\}^d$ and $z\in \Z^d$ we define $\xi^N(x,\eta)$ by
\begin{equation}
\xi^N(x,\eta):=\xi(y,\eta)
\end{equation}
It is well known (\cite{CapIof01})  that in the quenched regime (under the law $\P^{\xi^N(\eta)}_0$) as $\epsilon \downarrow 0$, $\epsilon X\big(t/\epsilon^2,\xi^N(\eta)\big)$ converges in law towards a Brownian Motion on $\Z^d$ with covariance matrix (effective diffusivity) $\sigma(\xi^N,\eta)$ (which is a random matrix on $X$, depending on the particular realization $\xi^N(\eta)$).

\subsubsection{The main theorem}
It is our purpose to prove the following theorem
\begin{Theorem}\label{dkhbdlhl8981nn2}
For $\mu$-almost all $\eta \in X$
\begin{equation}
\lim_{N\rightarrow +\infty}\sigma(\xi^N,\eta)=\sigma(\xi,\mu)
\end{equation}
\end{Theorem}
This result has already been given in \cite{CapIof01} when the jump rates are i.i.d. It is interesting to note that when the jump rates are i.i.d., D. Ioffe and P. Caputo have shown an exponential rate of convergence of $\sigma(\xi^N,\eta)$ towards $\E_{\mu}\big[\sigma(\xi^N,\eta)\big]$ as $N\rightarrow \infty$. 

\subsubsection{Proof}
We shall use the variational formula given in \cite{CapIof01}: for $l \in \R^d$
\begin{equation}\label{salksnlkn001}
^tl\sigma(\xi,\mu)l=\inf_{f \in L^2(\mu)}\sum_{i=1}^d \Big<\xi_i (l_i + D_i f)^2\Big>
\end{equation}
Let us also remind  the variational formula 
\begin{equation}\label{salksnliiikn001}
^tl\sigma(\xi^N,\eta)l=\inf_{f \in L^2(T^d_N)}|T^d_N|^{-1}\sum_{x\in T^d_N}\sum_{i=1}^d \xi_i^N(x,\eta) \Big(l_i + \nabla_i f(x)\Big)^2
\end{equation}
We will prove in the subsection  \ref{skasjsbjbd982} the following lemma which corresponds  to the variational formulation of $D^{-1}(\xi)$
\begin{Lemma}\label{suvsiuoswzvzvd82}
For $l\in \R^d$
\begin{equation}\label{dkjdkjbfbj8141}
^tl \sigma(\xi,\mu)^{-1}l=\inf_{H \in \S(X,\mu)}\sum_{i=1}^d \Big<\xi_i^{-1} (l_i + (\D H)_i)^2\Big>
\end{equation}
\end{Lemma}
It is also easy to prove that for $l\in \R^d$ (the proof is similar to the one of lemma \ref{dkjdkjbfbj8141})
\begin{equation}\label{shsjshbsjhb8001}
^tl \sigma(\xi^N,\eta)^{-1}l=\inf_{H\in \S(T^d_N)}N^{-d}\sum_{x\in T^d_N}\sum_{i=1}^d \xi^N_i(x,\eta)^{-1}\big(l_i+ (\Div H)_i \big)^2
\end{equation}
Let $l \in \R^d$.Let us apply theorem \ref{JHlkhbsljbj81} with $m=1$, $p=0$ and 
$\Phi(\eta, X^1)=\sum_{i=1}^d\xi_i(\eta) (l_i + X^1_i)^2$.
 By  the  Minkowski inequality and the uniform ellipticity condition \eref{djhdjhwhdheu1} one has that for $N\in \N^*$ and $\mu$-almost all $\eta \in X$, $v^1,v^2 \in \big(\L^2(T^d_N)\big)^2$
\begin{equation}
\big(\Psi_N(\eta,v^1)\big)^{\frac{1}{2}}-\big(\Psi_N(\eta,v^2)\big)^{\frac{1}{2}}\leq C \|v^1-v^2\|_{\L^2(T^d_N)}
\end{equation}
It follows that $\Phi$ is admissible and from the variational formulae
 \eref{salksnlkn001}, \eref{salksnliiikn001} and theorem \ref{JHlkhbsljbj81} one obtains that for $\mu$-almost all $\eta \in X$
\begin{equation}
\lim\sup_{N\rightarrow \infty}{^tl} \sigma(\xi^N,\eta)l\leq {^tl}\sigma(\xi,\mu)l
\end{equation}
Which gives the upper bound of theorem \ref{dkhbdlhl8981nn2}.
The proof of the lower bound is trivially similar: using Minkowski inequality, the uniform ellipticity condition \eref{djhdjhwhdheu1}, variational formulae
 \eref{dkjdkjbfbj8141}, \eref{shsjshbsjhb8001} and theorem \ref{JHlkhbsljbj81} one obtains that for $\mu$-almost all $\eta \in X$
\begin{equation}
\lim\sup_{N\rightarrow \infty}{^tl} \big(\sigma(\xi^N,\eta)\big)^{-1}l\leq {^tl} \big(\sigma(\xi,\mu)\big)^{-1}l
\end{equation}
Which gives the lower  bound of theorem \ref{dkhbdlhl8981nn2}.

\section{Proofs}
\subsection{Main results}\label{sdkjuhdjhdhh89112}
\subsubsection{Proof of  theorem \ref{suvsiuoswzvzvd81}}\label{skasjsbjbd981}
 It is trivial to check that $F^2_{pot}$, $F^2_{sol}$ and $\R^d$ are mutually orthogonal. Thus in order to  prove the Weyl decomposition \eref{suvsiuoswzvzvd81} it is sufficient to check that any element of $\L^2(X,\mu)$ orthogonal to $F^2_{pot}$ and $\R^d$ is an element of $F^2_{sol}$. Let $P$ be an element of $\L^2(X,\mu)$ orthogonal to $F^2_{pot}$ and $\R^d$. Since $P\perp F^2_{pot}$ it must verify
\begin{equation}
\sum_{i=1}^d D^*_i P=0
\end{equation}
By Lax-Milgram lemma for $m,n\in \{1,\ldots,d\}$, there exists $B^{m,n}\in F^2_{pot}$ such that
\begin{equation}\label{sskjsjdjd810001}
\sum_{i=1}^d D^*_i B^{m,n}_i=D_m^* P_n
\end{equation}
Let us define for $i,m,n\in \{1,\ldots,d\}$
\begin{equation}
H_{i,n,m}= B^{m,n}_i- B^{n,m}_i
\end{equation}
Let us define $Q\in \L^2(X,\mu)$ by for $n\in \{1,\ldots,d\}$
\begin{equation}
\begin{split}
Q_n=\sum_{i=1}^d H_{i,n,i}
\end{split}
\end{equation}
Since $B^{m,n}\in F^2_{pot}$ they can be approximated by gradient forms in $\L^2_{pot}$ and it is easy to deduce that $Q_n \in F^2_{sol}$. Moreover  for all $n\in \{1,\ldots,d\}$
\begin{equation}
\begin{split}
\sum_{k=1}^d D^*_k D_k Q_n=\sum_{k=1}^d \sum_{i=1}^d \Big(D^*_k D_k B^{i,n}_i- D^*_k D_k B^{n,i}_i\Big)
\end{split}
\end{equation}
Since $B^{m,n}\in F^2_{pot}$ it is easy to check by density that $\eta$-a.s., $D_k B^{i,n}_i=D_i B^{i,n}_k$ and $D_k B^{n,i}_i=D_i B^{n,i}_k$ thus from the equation \eref{sskjsjdjd810001} one obtains that $\eta$-a.s.
\begin{equation}
\begin{split}
\sum_{k=1}^d D^*_k D_k Q_n&= \sum_{i=1}^d D_i \sum_{k=1}^d \Big(D^*_k B^{i,n}_k(\eta)- D^*_k B^{n,i}_k(\eta)\Big)\\
&=\sum_{i=1}^d D_i \big(D^*_i P_n - D^*_n P_i\big)=\sum_{i=1}^d D_i^* D_i P_n - D_n \sum_{i=1}^d D^*_i P_i
\end{split}
\end{equation}
Using $\sum_{i=1}^d D^*_i P_i=0$ we obtain that $\eta$-a.s.
\begin{equation}
\begin{split}
\sum_{k=1}^d D^*_k D_k (Q_n- P_n )=0
\end{split}
\end{equation}
Combining this with $<Q_n-P_n>=0$, it follows by Lax Milgram lemma that $\eta-a.s.$, $Q_n=P_n$ and since $Q\in F^2_{sol}$ it follows that $P\in F^2_{sol}$ which concludes the proof of  theorem \ref{suvsiuoswzvzvd81}.

\subsubsection{Proof of equation \eref{dkdjdjdhjhd09033445b} of lemma  \ref{dddhdhdhhh8711b} in the discrete case}\label{skasjsbjbd984}
Let $p\in F^2_{sol}$. We will prove in this subsection that for $\mu$-almost all $\eta\in X$
\begin{equation}\label{dhkhhhhcbhh2333}
\lim_{N\rightarrow \infty} \|\Pi_N p(\eta)- \big(\Pi_N p(\eta)\big)_{sol}\|_{\L^2(T^d_N)}=0
\end{equation}
For $A\subset \R^d$ and $f\in L^2_{loc}(\Z^d)$ we will write
\begin{equation}\label{dkdjkdjdjkn1000222}
\|f\|_{L^2(A)}=\big(\sum_{x\in \Z^d\cap A} f(x)^2\big)^\frac{1}{2}
\end{equation}
Observe that to prove the equation \eref{dhkhhhhcbhh2333}, it is sufficient to prove the following lemma 
\begin{Lemma}\label{jedkjebdjh12333}
For $\mu$-almost all $\eta \in X$, there exists a sequence   $(K^N(x,\eta,M))_{M,N\in \N}$ of skew symmetric matrices with coefficients in $L^2(T^d_N)$ and a sequence of positive reals $h(M)$ such that $\lim_{M\rightarrow \infty } h(M)=0$ and for $M\geq 10$
\begin{equation}
\lim\sup_{N\rightarrow \infty} N^{-d/2}\|p( x,\eta)-\Div K^N(x,\eta,M)\|_{L^2([0,N(^d)} \leq h(M)
\end{equation}
\end{Lemma}
Let us now prove lemma \ref{jedkjebdjh12333}. Let $M\in \N,\;M\geq 10$, $M\leq 10^3 N$.
Since $p\in F^2_{sol}$, on obtains from Weyl decomposition \eref{suvsiuoswzvzvd80001} that for each $M$, there exists a $d\times d$  skew symmetric matrix $H^M$, with coefficients $H_{i,j}^M \in L^2(\mu)$, $(i,j)\in \{1,\ldots d\}^2$ such that $H_{i,j}^M=-H_{j,i}^M$ and
\begin{equation}\label{sdkndskdndsk20333}
\sum_{i=1}^d \big<|p- \Div H^M|^2\big> \leq 1/M^2
\end{equation}
and it is easy to check from the proof of \eref{suvsiuoswzvzvd80001} given in subsection \ref{skasjsbjbd981} that one can choose $H^M$ such that for all $(i,j)\in \{1,\ldots d\}^2$
\begin{equation}\label{sdkndskdndsk20333222333000}
 \big<| \nabla H^M_{i,j}|^2\big> \leq C_d \big<|p|^2\big> 
\end{equation}
Observe that by the ergodic theorem \ref{wlkkiudbu001}, $\eta$-a.s. 
\begin{equation}\label{sdkndskdndsk2333}
\lim\sup_{N\rightarrow \infty} \Big(N^{-d}\sum_{x\in \Z^d\cap [0,N(^d}|p( x,\eta)-\Div H^M(x,\eta)|^2\Big)^\frac{1}{2} \leq 1/M
\end{equation}
Let $g$ be a smooth increasing function on $\R$ such that $g(z)=1$ for $z\geq 1$ and $g(z)=0$ for $z\leq 1/2$ and let for $x\in [0,1]^d$
$\alpha_M(x)=g(M \operatorname{dist}(x,([0,1]^d)^c))$. Our candidate for $K^N$ will be the skew symmetric $T^d_N$-periodic matrix:
\begin{equation}
K^N(x,\eta,M)= \big(H^M(x,\eta)-[N/M]^{-d}\sum_{y\in [0,N/M(^d\cap \Z^d}H^M(y,\eta)\big) \alpha_M(x/N) \quad \text{on}\quad [0,N(^d\cap \Z^d
\end{equation}
Observe that  $\alpha_M(x)$ is null on an open neighborhood of $\R^d$ containing the boundary of $[0,1]^d$ and the coefficients of $K^N$ can be defined as  elements of $\L^2(T^d_N)$. 
Let us write
\begin{equation}\label{dwhbdlrebhb81333}
 J_1(N,M,\eta)=N^{-d/2}\|\Div H^M(x,\eta)-\Div K^N(x,\eta,M)\|_{L^2\big([0,N(^d\big)}
\end{equation}
Observe that $\eta$-a.s. ($\nabla \alpha_M$ standing for the discrete gradient of $\alpha_M$)
\[
\begin{split}
\Div K^N(x,\eta,M)=&\Div H^M(x,\eta) \alpha_M(x/N)\\&+ 
\big(H^M(x,\eta)-[N/M]^{-d}\sum_{y\in [0,N/M(^d\cap \Z^d}H^M(y,\eta)\big)  \nabla \big(\alpha_M(x/N)\big)
\end{split}
\]
 Thus  
\begin{equation}\label{sljhk7shjhfggl872333}
 J_1(N,M,\eta)\leq J_2(N,M,\eta)+ J_3(N,M,\eta)
\end{equation}
with
\begin{equation}
 J_2(N,M,\eta)=N^{-d/2}\|\Div H^M( x,\eta)(1-\alpha_M(x/N) )\|_{L^2\big([0,N(^d\big)}
\end{equation}
and
\begin{equation}
 J_3(N,M,\eta)=N^{-d/2}\|\big(H^M(x,\eta)-[N/M]^{-d}\sum_{y\in [0,N/M(^d\cap \Z^d}H^M(y,\eta)\big)  \nabla \big(\alpha_M(x/N)\big)\|_{L^2\big([0,N(^d\big)}
\end{equation}
Write 
\begin{equation}
A_M=\{x\in [0,N(^d\cap \Z^d \,:\, \min_{j\in \{1,\ldots d\}} \min(x_j,N-x_j)< N/M\}
\end{equation}
observe that
\begin{equation}
\begin{split}
 J_2(N,M,\eta)&\leq N^{-d/2}\|\Div H^M( x,\eta) \|_{L^2(A_M)}\\
&\leq C_d M^{-\frac{1}{2}} \big(\operatorname{Vol}(A_M)\big)^{-1/2}\|\Div H^M( x,\eta) \|_{L^2(A_M)}
\end{split}
\end{equation}
and by the ergodic theorem \ref{wlkkiudbu1}, $\eta$-a.s.\\ 
$ \big(\operatorname{Vol}(A_M)\big)^{-1/2}\|\Div H^M( x,\eta) \|_{L^2(A_M)}\rightarrow \big<(\Div H^M)^2\big>^\frac{1}{2}$ as $N\rightarrow \infty$.
 Thus $\eta$-a.s.
\begin{equation}\label{slweejhkhhsjhj7shl873333}
\lim \sup_{N\rightarrow \infty}  J_2(N,M,\eta) \leq C_d M^{-\frac{1}{2}} \big<p^2\big>^\frac{1}{2}
\end{equation}
Now let us prove that 
\begin{equation}\label{sdkndskdndsk20333222333}
\lim_{M \rightarrow \infty}\lim \sup_{N\rightarrow \infty}  J_3(N,M,\eta) =0
\end{equation}
Since $|\nabla (\alpha_M(x/N))|\leq C_d M/N$ one has 
\begin{equation}\label{kjsdhahol8771222333}
 J_3(N,M,\eta)\leq C_d N^{-1-d/2}M\sum_{m,n}\| H^M_{m,n}(x,\eta)-[N/M]^{-d}\sum_{y\in [0,N/M(^d\cap \Z^d}H^M_{m,n}(y,\eta)  \|_{L^2(A_M)} 
\end{equation} 
Let $I(M)=\{(i_1,\ldots,i_d)\in \{1,\ldots, M\}^d\,;\, \min_{j}\min(i_j-1,M-i_j)=0 \}$ and  write $\{B_i\}_{i\in I(M)}$ the set of cubes covering $A_M$ (the $N/M$-neighborhood of the border of $[0,N(^d$). More precisely for $i\in I(M)$,
\begin{equation}
B_i=\{x\in [0,N(^d\cap \Z^d\,:\, \max_j |x_j/N-(i_j-0.5)/M|\leq 1/(2M)\}
\end{equation}
By the equation \eref{kjsdhahol8771222333} one has
\begin{equation}\label{kjsdwehfaaadfdfahol8772333}
 J_3(N,M,\eta)^2\leq C_d M^2 \sum_{i \in I(M)}\sum_{m,n} K_i^{m,n}
\end{equation} 
with
\begin{equation}\label{kjsdhawedaaaffdhol8773333}
K_i^{m,n}=N^{-2-d}\big\|H^M_{m,n}(x,\eta)-[N/M]^{-d}\sum_{y\in [0,N/M(^d\cap \Z^d}H^M_{m,n}(y,\eta)\big\|^2_{L^2(B_i)}
\end{equation}
Now  using the inequality $(X+Y)^2\leq 2X^2+2Y^2$ observe that 
\begin{equation}\label{kjsdhfaaffaaffahol8774333}
\begin{split}
K_i^{m,n}\leq& 2 N^{-2-d} \big\|H^M_{m,n}(x,\eta)-\operatorname{Vol}(B_i))^{-1}\sum_{y\in B_i}H^M_{m,n}(y,\eta)\big\|^2_{L^2(B_i)}
\\& + 2 N^{-2} M^{-d} \Big([N/M]^{-d}\sum_{y\in [0,N/M(^d\cap \Z^d}H^M_{m,n}(y,\eta)-(\operatorname{Vol}(B_i))^{-1}\sum_{y\in B_i}H^M_{m,n}( y,\eta)\Big)^2
\end{split}
\end{equation}
By the Poincar\'{e} inequality one has 
\begin{equation}\label{kjsaaaadhahoaal8775333}
N^{-2-d} \big\|H^M_{m,n}(x,\eta)-\operatorname{Vol}(B_i))^{-1}\sum_{y\in B_i}H^M_{m,n}(y,\eta)\big\|^2_{L^2(B_i)} \leq C_d M^{-2} N^{-d}\big\|\nabla H^M_{m,n}( x,\eta)\big\|^2_{L^2(B_i)} 
\end{equation}
Thus
\begin{equation}\label{kjsdaaahahol87722333}
\begin{split}
 J_3(N,M,\eta)^2\leq& C_d M^{-1} \sum_{m,n}(\Vol(A_M))^{-1}\big\|\nabla H^M_{m,n}( x,\eta)\big\|^2_{L^2(A_M)}  \\+ C_d M^{2-d}\sum_{m,n} \sum_{i \in I(M)}& N^{-2} \Big([N/M]^{-d}\sum_{y\in [0,N/M(^d\cap \Z^d}H^M_{m,n}(y,\eta)-(\operatorname{Vol}(B_i))^{-1}\sum_{y\in B_i}H^M_{m,n}( y,\eta)\Big)^2
\end{split}
\end{equation} 
But by the ergodic theorem \ref{wlkkiudbu001}, $\eta-a.s.$ for all $i\in I(M)$
\begin{equation}\label{aslsldddssuxuc1222333}
\begin{split}
\lim_{N\rightarrow \infty} \Big([N/M]^{-d}\sum_{y\in [0,N/M(^d\cap \Z^d}H^M_{m,n}(y,\eta)&-(\operatorname{Vol}(B_i))^{-1}\sum_{y\in B_i}H^M_{m,n}( y,\eta)\Big)\\&= <H^M_{m,n}>-<H^M_{m,n}>=0
\end{split}
\end{equation}
It is important to observe that $H^M_{m,n}\in L^2(\mu)$ is sufficient to apply the ergodic theorem \ref{wlkkiudbu001} in order to obtain \eref{aslsldddssuxuc1222333}. In the continuous we will not consider an approximation of $p$ but its direct primitive which is not in $L^2(\mu)$ explaining why the ergodic theorem will not be applied directly.\\
It follows that (using \eref{sdkndskdndsk20333222333000})
\begin{equation}\label{kjsdhahol87723222333}
\begin{split}
 \lim \sup_{N\rightarrow \infty }J_3(N,M,\eta)^2&\leq C_d M^{-1}\sum_{m,n} <(\nabla H^M_{m,n})^2> 
\\&\leq C_d M^{-1} <p^2> 
\end{split}
\end{equation} 
And taking the limit $M\rightarrow \infty$ one obtains the equation \eref{sdkndskdndsk20333222333}.
Now combining equations \eref{sdkndskdndsk20333}, \eref{dwhbdlrebhb81333}, \eref{sljhk7shjhfggl872333}, \eref{slweejhkhhsjhj7shl873333} and \eref{sdkndskdndsk20333222333}   one obtains lemma \ref{jedkjebdjh1222}.

\subsubsection{Proof of equation \eref{dkdjdjdhjhd090334b} of lemma  \ref{dddhdhdhhh8711b} in the discrete case}\label{skasjsbjbd983}
The proof of  equation \eref{dkdjdjdhjhd090334b} being similar to the one of equation \eref{dkdjdjdhjhd09033445b} we will just give its idea.
Let $v\in F^2_{pot}$. We have to prove that for $\mu$-almost all $\eta\in X$
\begin{equation}\label{dvbdvdjdvjdvh89918000}
\lim_{N\rightarrow \infty} \|\Pi_N v(\eta)- \big(\Pi_N v(\eta)\big)_{pot}\|_{\L^2(T^d_N)}=0
\end{equation}
Observe that to prove the equation \eref{dvbdvdjdvjdvh89918000} it is sufficient to prove the following lemma 
\begin{Lemma}\label{jedkjebdjh1222}
For $\mu$-almost all $\eta \in X$, there exists a sequence $(G^N(x,\eta,M))_{M,N\in \N}$ of functions in $\L^2(T^d_N)$  and sequence of positive reals $h(M)$ such that $\lim_{M\rightarrow \infty } h(M)=0$ and for $M\geq 10$
\begin{equation}
\lim\sup_{N\rightarrow \infty} \Big(N^{-d}\sum_{x\in \Z^d\cap [0,N(^d}|v( x,\eta)-\nabla G^N(x,\eta,M)|^2\Big)^\frac{1}{2} \leq h(M)
\end{equation}
\end{Lemma}
Let us now prove lemma \ref{jedkjebdjh1222}. Let $M\in \N,\;M\geq 10$, $M\leq 10^3 N$. Since $v\in F^2_{pot}$,  there exists $u^M\in L^2(\mu)$ such that
\begin{equation}\label{sdkndskdndsk20}
\sum_{i=1}^d \big<|v_i-D_i u^M|^2\big> \leq 1/M^2
\end{equation}
Observe that by the ergodic theorem \ref{wlkkiudbu001}, $\eta$-a.s. 
\begin{equation}\label{sdkndskdndsk2}
\lim\sup_{N\rightarrow \infty} \Big(N^{-d}\sum_{x\in \Z^d\cap [0,N(^d}|v( x,\eta)-\nabla u^M(x,\eta)|^2\Big)^\frac{1}{2} \leq 1/M
\end{equation}
Defining $\alpha_M(x)$ as in the subsection \ref{skasjsbjbd984} our candidate
 for $G^N$ will be the $\L^2(T^d_N)$ periodic function with value ($[N/M]$ being the integer part of $N/M$)
\begin{equation}
G^N(x,\eta,M)= \big(u^M(x,\eta)-[N/M]^{-d}\sum_{y\in [0,N/M(^d\cap \Z^d}u^M(y,\eta)\big) \alpha_M(x/N) \quad \text{on}\quad [0,N(^d\cap \Z^d
\end{equation}
From this point the proof of lemma \ref{jedkjebdjh1222} is trivially similar to the one given in subsection \ref{skasjsbjbd984}.

\subsubsection{Proof of equation \eref{dkdjdjdhjhd090334b} of lemma  \ref{dddhdhdhhh8711b} in the continuous case}\label{eicnboiecbb81}
Let $v\in F^2_{pot}$. We will prove in this subsection that for $\mu$-almost all $\eta\in X$
\begin{equation}\label{djhdjheuu1}
\lim_{N\rightarrow \infty} \|\Pi_N v(\eta)- \big(\Pi_N v(\eta)\big)_{pot}\|_{\L^2(T^d_N)}=0
\end{equation}
Observe that to prove the equation \eref{djhdjheuu1} it is sufficient to prove the following lemma 
\begin{Lemma}\label{jedkjebdjh1}
For $\mu$-almost all $\eta \in X$, there exists a sequence $(G^N(x,\eta,M))_{M,N\in \N}$ of functions in $H^1(T^d_1)$  and sequence of positive reals $h(M)$ such that $\lim_{M\rightarrow \infty } h(M)=0$ and for $M\geq 10$
\begin{equation}
\lim\sup_{N\rightarrow \infty} \|v(N x,\eta)-\nabla G^N(x,\eta,M)\|_{L^2([0,1]^d)} \leq h(M)
\end{equation}
\end{Lemma}
Let us now prove lemma \ref{jedkjebdjh1}. Since $v\in F^2_{pot}$, for almost all $\eta$, $v$ admits the following representation $v(x,\eta)=\nabla_x u(x,\eta)$, where $u(x,\eta)$ is an element of $H^1_{loc}(\R^d)$ (see subsection \ref{jdhdbbhdhii101}).\\
Let $M\in \N,\;M\geq 10$.
Let $z\rightarrow g(z)$ be a smooth increasing function on $\R$ such that $g=1$ for $z\geq 1$ and $g=0$ for $z\leq 1/2$ and let for $x\in [0,1]^d$
$\alpha_M(x)=g\big(M \operatorname{dist}(x,([0,1]^d)^c)\big)$. Our candidate for $G^N$ will be the $H^1(T^d_1)$ periodic function with value
\begin{equation}
G^N(x,\eta,M)= N^{-1}\big(u(Nx,\eta)-M^d\int_{[0,1/M]^d}u(Ny,\eta)dy\big) \alpha_M(x) \quad \text{on}\quad [0,1]^d
\end{equation}
Observe that since $\alpha_M(x)$ is null on an open neighborhood of $\R^d$ containing the boundary of $[0,1]^d$, $G^N$ can be defined as an element of $H^1(T^d_1)$. Let us write
\begin{equation}\label{dkdjkdjdjkn1}
 J_1(N,M,\eta)=\|\nabla u(N x,\eta)-\nabla G^N(x,\eta,M)\|_{L^2([0,1]^d)}
\end{equation}
Observe that $$\nabla G^N(x,\eta,M)=\nabla u(Nx,\eta) \alpha_M(x)+ N^{-1}\big(u(Nx,\eta)-M^d\int_{[0,1/M]^d}u(Ny,\eta)dy\big)  \nabla \alpha_M(x)$$ 
Thus 
\begin{equation}\label{sljhk7shl872}
 J_1(N,M,\eta)\leq J_2(N,M,\eta)+ J_3(N,M,\eta)
\end{equation}
with
\begin{equation}
 J_2(N,M,\eta)=\|\nabla u(N x,\eta) (1-\alpha_M(x) )\|_{L^2([0,1]^d)} 
\end{equation}
and
\begin{equation}
 J_3(N,M,\eta)=N^{-1}\| \big(u(Nx,\eta)-M^d\int_{[0,1/M]^d}u(Ny,\eta)dy\big)  \nabla \alpha_M(x) \|_{L^2([0,1]^d)} 
\end{equation}
Write 
\begin{equation}
A_M=\{x\in [0,1]^d\,:\, \min_{j\in \{1,\ldots d\}} \min(x_j,1-x_j)\leq 1/M\}
\end{equation}
Observe that
\begin{equation}
\begin{split}
 J_2(N,M,\eta)&\leq \|\nabla u(N x,\eta) \|_{L^2(A_M)}\\
&\leq C_d M^{-\frac{1}{2}} \big(\operatorname{Vol}(A_M)\big)^{-1/2}\|\nabla u(N x,\eta) \|_{L^2(A_M)}
\end{split}
\end{equation}
and by the ergodic theorem \ref{wlkkiudbu1}, $\eta$-a.s. $ \big(\operatorname{Vol}(A_M)\big)^{-1/2}\|\nabla u(N x,\eta) \|_{L^2(A_M)}\rightarrow \big<v^2\big>^\frac{1}{2}$ as $N\rightarrow \infty$. Thus $\eta$-a.s.
\begin{equation}\label{sljhk7shl873}
\lim \sup_{N\rightarrow \infty}  J_2(N,M,\eta) \leq C_d M^{-\frac{1}{2}} \big<v^2\big>^\frac{1}{2}
\end{equation}
Now let us prove that $\eta-a.s.$
\begin{equation}\label{salkjc2b781}
\lim_{M \rightarrow \infty}\lim \sup_{N\rightarrow \infty}  J_3(N,M,\eta) =0
\end{equation}
Since $|\nabla \alpha_M|\leq C_d M$ one has 
\begin{equation}\label{kjsdhahol8771}
 J_3(N,M,\eta)\leq C_d N^{-1} \|M \big(u(Nx,\eta)-M^d\int_{[0,1/M]^d}u(Ny,\eta)dy\big)  \|_{L^2(A_M)} 
\end{equation} 
Let $I(M)=\{(i_1,\ldots,i_d)\in \{1,\ldots, M\}^d\,;\, \min_{j}\min(i_j-1,M-i_j)=0 \}$ and $\{B_i\}_{i\in I(M)}$ the set of cubes covering $A_M$ (the $1/M$-neighborhood of the border of $[0,1]^d$). More precisely for $i\in I(M)$,
\begin{equation}\label{sjbsbdcjbdfufb881112}
B_i=\{x\in [0,1]^d\,:\, \max_j |x_j-(i_j-0.5)/M|\leq 1/(2M)\}
\end{equation}
By the equation \eref{kjsdhahol8771} one has
\begin{equation}\label{kjsdhahol8772}
 J_3(N,M,\eta)^2\leq C_d M^2 \sum_{i \in I(M)} K_i
\end{equation} 
with
\begin{equation}\label{kjsdhahol8773}
K_i=N^{-2}\int_{B_i} \big(u(Nx,\eta)-M^d\int_{[0,1/M]^d}u(Ny,\eta)dy\big)^2\, dx 
\end{equation}
Now using the inequality $(X+Y)^2\leq 2X^2+2Y^2$ observe that 
\begin{equation}\label{kjsdhahol8774}
\begin{split}
K_i\leq& 2 N^{-2} \int_{B_i} \big(u(N x,\eta)-(\operatorname{Vol}(B_i))^{-1}\int_{B_i}u(N y,\eta)dy\big)^2\, dx\\& + 2 N^{-2}(\operatorname{Vol}(B_i))^{-1} \Big(\int_{[0,1/M]^d}u(N y,\eta)dy-\int_{B_i}u(N y,\eta)dy\Big)^2
\end{split}
\end{equation}
By the Poincar\'{e} inequality one has 
\begin{equation}\label{kjsdhahol8775}
N^{-2} \int_{B_i} \big(u(N x,\eta)-(\operatorname{Vol}(B_i))^{-1}\int_{B_i}u(N y,\eta)dy\big)^2\, dx \leq C_d M^{-2} \int_{B_i} \big(v(N x,\eta)\big)^2 dx
\end{equation}
Thus
\begin{equation}\label{kjsdhahol87722}
\begin{split}
 J_3(N,M,\eta)^2\leq& C_d M^{-1} (\Vol(A_M))^{-1}\int_{A_M} \big(v(N x,\eta)\big)^2 dx \\&+ C_d M^2 \sum_{i \in I(M)} (\operatorname{Vol}(B_i))^{-1}N^{-2} \Big(\int_{[0,1/M]^d}u(N y,\eta)dy-\int_{B_i}u(N y,\eta)dy\Big)^2
\end{split}
\end{equation} 
It shall be proven in the paragraph \ref{aslslssuxuc1} that by the
 ergodic theorem $\eta-a.s.$ for all $i\in I(M)$
\begin{equation}\label{aslsldddssuxuc1}
\lim_{N\rightarrow \infty} N^{-1}\big|\int_{[0,1/M]^d}u(N y,\eta)dy-\int_{B_i}u(N y,\eta)dy\big|=0
\end{equation}
It follows that $\eta-a.s.$
\begin{equation}\label{kjsdhahol87723}
\begin{split}
 \lim \sup_{N\rightarrow \infty }J_3(N,M,\eta)^2\leq C_d M^{-1} <v^2> 
\end{split}
\end{equation} 
And taking the limit $M\rightarrow \infty$ one obtains the equation \eref{salkjc2b781}.
Now combining equations \eref{dkdjkdjdjkn1}, \eref{sljhk7shl872}, \eref{sljhk7shl873}  one obtains  lemma \ref{jedkjebdjh1} with $h(M)=C_d M^{-1/2} <v^2>^{1/2} $.

\paragraph{Proof of equation \eref{aslsldddssuxuc1}}\label{aslslssuxuc1}
Let $a, b \in I(M)\times I(M)$, $q\in \N^{*}$ and $w\in \{1,\ldots,d\}$ such that $b=a+q e_w$. Observing that 
$$u(N (x+q  e_w/M),\eta)-u(N x,\eta)=\int_{t\in [0,1]}v(N (x+q t e_w/M) ,\eta).(q e_w/M)\,dt
$$
we obtain
\begin{equation}\label{shshsvdfchvfvf1}
\begin{split}
N^{-1}\int_{B_b}u(N x,\eta)dx=&N^{-1} \int_{B_a}u(N x,\eta)dx\\&+\int_{t\in [0,1]}\int_{B_a}v(N (x+q t e_w/M ),\eta).(q  e_w/M)dx\,dt
\end{split}
\end{equation}
Now let us write $\partial^w B_a$ the lower face of the cube $B_a$ orthogonal to $e_w$: $$\partial^w B_a:=\{x\in B_a\,:\, x.e_w=(a_w-0.5)/M-1/(2M)\}$$
Now we decompose $x\in B_a$ as $x= x^w+y_w e_w/M$ with $x^w\in \partial^w B_a$ and $y_w \in [0,1]$, using this change of variable we obtain 
\begin{equation}\label{gdgdghsdgds11}
\begin{split}
\int_{t\in [0,1]}&\int_{B_a}v(N (x+q t e_w/M) ,\eta).(q  e_w/M)dx\,dt\\&=\int_{t\in [0,1]}\int_{y_w \in [0,1], x^w \in \partial^w B_a}v(N (x^w+(tq+y_w) e_w/M ),\eta). (qe_w/M^2)\, dx^w\,dt\,dy_w
\end{split}
\end{equation}
Now using the change of variable $s=(tq+y_w) /M$ we obtain from \eref{gdgdghsdgds11} that
\begin{equation}\label{sgsgsvsv88112}
\begin{split}
\int_{t\in [0,1]}&\int_{B_a}v(N (x+q t e_w/M) ,\eta).(q  e_w/M)dx\,dt\\&= \int_{y_w\in [0,1]}\int_{s\in [y_w/M,(q+y_w)/M]}\int_{ x^w \in \partial^w B_a}v(N (x^w+s e_w) ,\eta). (e_w/M)\, dx^w\,ds\,dy_w
\\&= \int_{s\in [1/M,q/M]}\int_{ x^w \in \partial^w B_a}v(N (x^w+s e_w ),\eta). (e_w/M)\, dx^w\,ds
\\&+\int_{s\in [0,1/M]}\int_{ x^w \in \partial^w B_a}v(N (x^w+s e_w ),\eta). e_w s\, dx^w\,ds
\\&+\int_{s\in [q/M,(q+1)/M]}\int_{ x^w \in \partial^w B_a}v(N (x^w+s e_w) ,\eta). e_w (1+q-Ms)\, dx^w\,ds
\end{split}
\end{equation}
Thus using the decomposition $[1/M,q/M]=\cup_{k=1}^{q-1} \big[k/M,(k+1)/M\big]$, it follows from \eref{sgsgsvsv88112} and \eref{shshsvdfchvfvf1}  by obvious change of variables that
\begin{equation}\label{hddkhfhhhf88812221}
\begin{split}
N^{-1}\int_{B_b}u(N x,\eta)dx
=&N^{-1} \int_{B_a}u(N x,\eta)dx+M^{-1}\sum_{k=1}^{q-1}\int_{B_a+k e_w}v(N x,\eta).e_w \,dx\\
&+ \int_{B_a}(x.e_w-a.e_w/M+1/M)v(N x ,\eta).e_w \,dx\\
&+ \int_{B_b}(-x.e_w+b.e_w/M)v(N x ,\eta).e_w \,dx
\end{split}
\end{equation}
Since $<v>=0$ one has $\eta$-a.s.
\begin{equation}\label{hddkhfhhhf88812222}
\begin{split}
\lim_{N\rightarrow \infty} \sum_{k=1}^{q-1}\int_{B_a+k e_w}v(N x,\eta).e_w \,dx=0
\end{split}
\end{equation}
Now   for $P\in \N^{*}$ ($P\geq M^2$) write $E_k=\{x\in B_b\,:\,(k-1)/P<-x.e_w+b.e_w/M \leq k/P \}$. Note that $(E_k)_{1\leq k\leq P}$ is a partition of $B_b$ (observe that by equation \eref{sjbsbdcjbdfufb881112} $b/M$ is not the center of $B_b$ but the upper edge of the cube). Thus we obtain
\begin{equation}
\begin{split}
 \int_{B_b}(-x.e_w+b.e_w/M)v(N x ,\eta).e_w \,dx\leq &\sum_{k=1}^P \int_{E_k}(k/P)v(N x ,\eta).e_w \,dx
\\&+ 1/P   \int_{B_b}|v(N x ,\eta).e_w| \,dx
\end{split}
\end{equation}
It follows by the ergodic theorem that $\eta$-a.s.
\begin{equation}
\begin{split}
\lim\sup_{N\rightarrow \infty} |\int_{B_b}(-x.e_w+b.e_w/M)v(N x ,\eta).e_w \,dx| \leq  P^{-1} M^{-d} <|v.e_w|>
\end{split}
\end{equation}
And taking the limit $P\rightarrow \infty$ one obtains that
\begin{equation}\label{hddkhfhhhf88812223}
\begin{split}
\lim\sup_{N\rightarrow \infty} |\int_{B_b}(-x.e_w+b.e_w/M)v(N x ,\eta).e_w \,dx| =0
\end{split}
\end{equation}
Similarly one obtains that 
\begin{equation}\label{hddkhfhhhf88812224}
\begin{split}
\lim\sup_{N\rightarrow \infty} |\int_{B_a}(x.e_w-a.e_w/M+1/M)v(N x ,\eta).e_w \,dx| =0
\end{split}
\end{equation}
From \eref{hddkhfhhhf88812221}, \eref{hddkhfhhhf88812222}, \eref{hddkhfhhhf88812223} and \eref{hddkhfhhhf88812224}  one deduces that for $a, b \in I(M)\times I(M)$, $q\in \N^{*}$ and $w\in \{1,\ldots,d\}$ such that $b=a+q e_w$. one has
\begin{equation}
\begin{split}
\lim\sup_{N\rightarrow \infty} |N^{-1}\int_{B_b}u(N x,\eta)dx-N^{-1} \int_{B_a}u(N x,\eta)dx|=0
\end{split}
\end{equation}
And since any two distinct points of $I(M)$ can be connected by a finite number of steps of such translations one obtains the equation \eref{aslsldddssuxuc1}

\subsubsection{Proof of equation \eref{dkdjdjdhjhd09033445b} of lemma  \ref{dddhdhdhhh8711b} in the continuous case}\label{eicnboiecbb82}
Let $p\in F^2_{sol}$. We will prove in this subsection that for $\mu$-almost all $\eta\in X$
\begin{equation}\label{dhkhhhhcbhh2}
\lim_{N\rightarrow \infty} \|\Pi_N p(\eta)- \big(\Pi_N p(\eta)\big)_{sol}\|_{\L^2(T^d_N)}=0
\end{equation}
Observe that to prove the equation \eref{dhkhhhhcbhh2}, it is sufficient to prove the following lemma
\begin{Lemma}\label{jedkjebdjh12}
For $\mu$-almost all $\eta \in X$, there exists a sequence   $(K^N(x,\eta,M))_{M,N\in \N}$ of skew symmetric matrices with coefficients in $H^1(T^d_1)$ and a sequence of positive reals $h(M)$ such that $\lim_{M\rightarrow \infty } h(M)=0$ and for $M\geq 10$
\begin{equation}
\lim\sup_{N\rightarrow \infty} \|p(N x,\eta)-\Div K^N(x,\eta,M)\|_{L^2([0,1]^d)} \leq h(M)
\end{equation}
\end{Lemma}
Since $p\in F^2_{sol}$, it is easy to prove from Weyl decomposition that there exists a finite sequence $h_{i,j} \in F^2_{pot}$, $(i,j)\in \{1,\ldots d\}^2$ such that $h_{i,j}=-h_{j,i}$ and $\eta$-a.s., $(p)_i=\sum_{j=1}^d h_{i,j}.e_j$ ($F_{pot}$ is orthogonal to the set of such vectors and any element of $\L^2(X,\mu)$ orthogonal to the set of such vectors is in $\L^2_{pot}(X,\mu)$).
Write $H_{i,j}$ the scalar potentials associated to $h_{i,j}$, then it follows that $\eta$-a.s., $\Div_x H_{i,j}(x,\eta)=h_{i,j}(x,\eta)$. Thus $H$ is a $d\times d$ skew symmetric matrix with elements in $H^1_{loc}(\R^d)$ such that $\eta$-a.s.,  $p(x,\eta)=\Div H(x,\eta)$.\\
For $M\in \N$, $M\geq 10$, defining $\alpha_M(x)$ as in the subsection \ref{eicnboiecbb81} our candidate for $K^N$ will be the skew symmetric $T^d_1$-periodic matrix:
\begin{equation}
K^N(x,\eta,M):= N^{-1}\big(H(Nx,\eta)-M^d\int_{[0,1/M]^d}H(Ny,\eta)dy\big) \alpha_M(x) \quad \text{on}\quad [0,1]^d
\end{equation}
From this point the proof of lemma \ref{jedkjebdjh12} is trivially similar to the one given in subsection \ref{eicnboiecbb81}.

\subsection{Applications}
In this subsection we will prove theorem \ref{dkhdbldbdhjb81}. We will first prove equation \eref{dsdjhdldhi8101}, the equations \eref{ddlkjdj1} and \eref{HoCoeDNorPdddre12} will be implied by the first one.
\subsubsection{Proof of  the variational formula \eref{dsdjhdldhi8101}}\label{hdodkdbbii871}
Let us write $v_l^{E}\in F^2_pot$ the solution of the equation \eref{wlbuasssbu1} and (using the linearity of $v_l^E$ in $l$ one can define $v_.^{E}$ as a matrix by $v_.^E l=v_l^E$). 
Let us first prove the following lemma
\begin{Lemma}\label{ksdjhbfclkb981}
\begin{equation}
\sigma(a,-E,\mu)={^t \sigma(a,E,\mu)}
\end{equation}
\end{Lemma}
\begin{proof}
We shall adapt the proof given by J. R. Norris \cite{Nor97} for the periodic case. 
Let $l,k\in \R^d$. Since $(a+E)(l+v^{E}_l)\in \L^2_{sol}(X,\mu)$ and $(a-E)(k+v^{-E}_k)\in \L^2_{sol}(X,\mu)$, by the Weyl decomposition \eref{aslkjhbiqsuu71},  there exists $q,h\in F^2_{sol}$ and $t,s \in \R^d$ such that 
\begin{equation}
t-q=(a+E)(l+v^E_l)
\end{equation}
and
\begin{equation}
s-h=(a-E)(k+v^{-E}_k)
\end{equation}
Observe that by integration with respect to the measure $\mu$, one obtains that
\begin{equation}
t=\sigma(a,E,\mu)l
\end{equation}
and
\begin{equation}
s=\sigma(a,-E,\mu)k
\end{equation}
For $f,g \in \L^2(X,\mu)$ we write $\Big<f,g\Big>=\Big<{^t f}g\Big>$. Then observe that
\begin{equation}
\begin{split}
\Big< \sigma(a,-E,\mu) k, l\Big>&=\Big< s, l\Big>=\Big< s-h, l+v_l^E\Big>
\\&=\Big< (a-E)(k+v^{-E}_k), l+v_l^E\Big>=\Big< k+v^{-E}_k,(a+E) (l+v_l^E)\Big>
\\&=\Big< k+v^{-E}_k,t-q\Big>=\Big< k,t\Big>
\\&=\Big<  k,\sigma(a,E,\mu) l\Big>
\end{split}
\end{equation}
Which proves that $^t\sigma(a,-E,\mu)=\sigma(a,E,\mu)$ and henceforth the lemma.
\end{proof}
Let $\xi,l \in \R^d$, we will now prove that 
\begin{equation}\label{dsdjhdldhi81010}
\begin{split}
|\xi-\sigma(a,E,\mu) l|_{\sigma_{\sym}^{-1}(a,E,\mu)}^2=\inf_{\psi,p\in F_{pot}\times F_{sol}}\Big<|\xi-p-(a+E)(l-\psi)|_{a^{-1}}^2\Big>
\end{split}
\end{equation}
We will write $\sigma_{\sym}$ is the symmetric part of $\sigma(a,E,\mu)$. Let us define
\begin{equation}
\psi_0:=v^E_.\big(l+\frac{1}{2} \sigma_{\sym}^{-1} (\xi-\sigma l))\big)-v^{-E}_. \frac{1}{2}\sigma_{\sym}^{-1} (\xi-\sigma l)
\end{equation}
and
\begin{equation}
p_{0}:=\xi- (a+E) (l-\psi_0)
-a \big(I_d- v^{-E}_.\big)\sigma_{\sym}^{-1} (\xi-\sigma l)
\end{equation}
Observe also that since 
\begin{equation}\label{dhkjdhdbhdd1b81}
\xi-p_0-(a+E)(l-\psi_0)= a (I_d- v^{-E}_.)\sigma_{\sym}^{-1}(\xi-\sigma l)
\end{equation}
And (using lemma \ref{ksdjhbfclkb981})
\begin{equation}
\big<{^t(}I_d- v^{-E}_.)a(I_d- v^{-E}_.)\big>=\sigma_{\sym}(a,-E,\mu)=\sigma_{\sym}(a,E,\mu)=\sigma_{\sym}
\end{equation}
One obtains that
\begin{equation}\label{dsdjhdldhi2281010}
\begin{split}
\Big<|\xi-p_0-(a+E)(l-\psi_0)|_{a^{-1}}^2\Big>=|\xi-\sigma(a,E,\mu) l|_{\sigma_{\sym}^{-1}(a,E,\mu)}^2
\end{split}
\end{equation}
Moreover $\psi_0 \in F^2_{pot}$ and $p_{0}\in F^2_{sol}$ since 
\begin{equation}
\begin{split}
p_{0}=&\xi-(a+E)l-a\sigma_{\sym}^{-1} (\xi-\sigma l)
\\&+ (a+E)(I_d+v_.^E)\big(l+\frac{1}{2} \sigma_{\sym}^{-1} (\xi-\sigma l))  \\&
+(a-E)(I_d+v_{.}^{-E}) \frac{1}{2} \sigma_{\sym}^{-1} (\xi-\sigma l)
\end{split}
\end{equation}
And by the equation $\eref{dhkjdhdbhdd1b81}$,  $a^{-1}\big(\xi-p_0-(a+E)(l-\psi_0)\big)$ is orthogonal in $\L^2(X,\mu)$ to 
$F^2_{sol}$ and the space $\{(a+E)v\,:\, v\in F^2_{pot} \}$, it follows that the variational formula \eref{dsdjhdldhi81010} is valid and the minimum is reached at $p_0$ and $\psi_0$.

\subsubsection{Proof of the variational formulas \eref{ddlkjdj1} and \eref{HoCoeDNorPdddre12}}
One obtains the variational formula \eref{ddlkjdj1} from the variational formula \eref{dsdjhdldhi8101} by observing that
\begin{equation}
\inf_{\xi\in \R^d, \xi\perp l}|\xi-\sigma(a,E,\mu) l|_{\sigma_{\sym}^{-1}(a,E,\mu)}^2={^t l}\sigma_{\sym}(a,E,\mu)l
\end{equation}
One obtains the variational formula \eref{HoCoeDNorPdddre12} by taking $l=0$ in \eref{dsdjhdldhi8101}.

\subsubsection{Proof of  lemma \ref{suvsiuoswzvzvd82}}\label{skasjsbjbd982}
Gift $\L^2(X,\mu)$ with the scalar product $(f,g)_H=\sum_{i=1}^d <\xi_i f_i g_i>$ to obtain an Hilbert space. By the variational formula \eref{salksnlkn001}, $^tlD(\xi,\mu)l$ is the norm of the H-orthogonal projection of $l$ on the subspace of $\L^2(X,\mu)$ H-orthogonal to $F^2_{pot}$. It follows that there exists an unique $v_l \in F^2_{pot}$ linear in $l$ realizing the minimum of  \eref{salksnlkn001} and such that $(l+v_l)$ is H-orthogonal to $F^2_{pot}$.\\
Thus the vector for $l\in \R^d$ the vector field $p_l$ defined by
\begin{equation}
p_l=\xi_i (I_d + v_.)\big(D(\xi,\mu)\big)^{-1}l -l
\end{equation}
verifies $<p_l>=0$ and is orthogonal to $F^2_{pot}$, thus by the theorem \ref{suvsiuoswzvzvd81}, it is an element of $F^2_{sol}$. Moreover observing that 
\begin{equation}\label{dk444jdkjbfbj8141}
 \Big<\xi_i^{-1} (l_i + (p_l)_i)^2\Big>= {^tl D(\xi,\mu)^{-1}l}
\end{equation}
and since the vector $q$ defined by $q_i=\xi_i^{-1} (l_i + (p_l)_i)=(I_d + v_.)\big(D(\xi,\mu)\big)^{-1}l$ is orthogonal to $F^2_{sol}$ one obtains that the variational formula \eref{suvsiuoswzvzvd82} is true and that its minimum is reached at $p_l$.

\paragraph*{Acknowledgments}
Part of this work was supported by the Aly Kaufman fellowship. The author would like to thank Dmitry Ioffe for his hospitality during his stay at the Technion, for suggesting this problem and for stimulating and helpful discussions. Thanks are also due to the referee for many useful comments.


\begin{thebibliography}{MFGW89}

\bibitem[BLP78]{BeLiPa78}
A.~Bensoussan, J.~L. Lions, and G.~Papanicolaou.
\newblock {\em Asymptotic analysis for periodic structure}.
\newblock North Holland, Amsterdam, 1978.

\bibitem[CI01]{CapIof01}
Pietro Caputo and Dmitry Ioffe.
\newblock Finite volume approxiamtion of the effective diffusion matrix: the
  case of independent bond disorder.
\newblock {\em preprint ar{X}iv:math.{PR}/0110215}, 2001.

\bibitem[DGI00]{DGI00}
J-D. Deuschel, G.~Giacomin, and D.~Ioffe.
\newblock Large deviations and concentration properties for $\nabla \phi$
  interface models.
\newblock {\em Probab. Theory Related Fields}, 117:49--111, 2000.

\bibitem[DS67]{DuScha67}
N.~Dunford and J.~T. Schwartz.
\newblock {\em Linear {O}perators, {P}art {I}}.
\newblock Interscience plublishers, 1967.

\bibitem[FP94]{FaPa94}
A.~Fannjiang and G.C. Papanicolaou.
\newblock Convection enhanced diffusion for periodic flows.
\newblock {\em SIAM J. Appl. Math.}, 54:333--408, 1994.

\bibitem[FP96]{FaPa96}
A.~Fannjiang and G.C. Papanicolaou.
\newblock Diffusion in turbulence.
\newblock {\em Probab. Theory Related Fields}, 105(3):279--334, 1996.

\bibitem[GOS01]{GOS01}
G.B. Giacomin, S.~Olla, and H.~Spohn.
\newblock Equilibrium fluctuation for a {G}inzburg-{L}andau $\nabla \phi$
  interface model.
\newblock {\em Annals of Probability}, 2001.

\bibitem[JKO91]{JiKoOl91}
V.~V. Jikov, S.~M. Kozlov, and O.~A. Oleinik.
\newblock {\em Homogenization of Differential Operators and Integral
  Functionals}.
\newblock Springer-Verlag, 1991.

\bibitem[Koz80]{Ko80}
S.M. Kozlov.
\newblock Averaging of random operators.
\newblock {\em Math USSR Sbornik}, 37:167--180, 1980.

\bibitem[Koz85]{Ko85}
S.M. Kozlov.
\newblock The method of averaging and walks in inhomogeneous environments.
\newblock {\em Russian Math. Surveys}, 2(40):73--145, 1985.

\bibitem[KV86]{KiVa86}
C.~Kipnis and S.R.S. Varadhan.
\newblock Central limit theorem for additive functional of reversible markov
  processes and application to simple exclusion.
\newblock {\em Comm. Math. Phys.}, 104:1--19, 1986.

\bibitem[LOV01]{LOV00}
C.~Landim, S.~Olla, and S.R.S. Varadhan.
\newblock Finite-dimensional approximation of the self-diffusion coefficient.
\newblock {\em Comm. Math. Phys.}, 224:302--321, 2001.

\bibitem[MFGW89]{MaFe89}
A.~De Masi, P.~A. Ferrari, S.~Goldstein, and W.D. Wick.
\newblock An invariance principle for reversible markov processes. application
  to random motions in random environments.
\newblock {\em Journal of Statistical Physics}, 55(3/4):787--855, 1989.

\bibitem[Nor97]{Nor97}
J.R. Norris.
\newblock Long-time behaviour of heat flow: Global estimates and exact
  asymptotics.
\newblock {\em Arch. Rational Mech. Anal.}, 140:161--195, 1997.

\bibitem[Oll94]{Ol94}
S.~Olla.
\newblock {\em Homogenization of Diffusion Processes in Random Fields}.
\newblock Ecole Polytechnique, 1994.
\newblock Cours Ecole Polytechnique.

\bibitem[OS95]{OsSa95}
H.~Osada and T.~Satoh.
\newblock An invariance principle for non-symmetric markov processes and
  reflecting diffusions in random domains.
\newblock {\em Probab. Theory Relat. Fields}, 101:45--63, 1995.

\bibitem[Pia02]{Piat02}
Andrey Piatnitski.
\newblock {\em Private discussion}.
\newblock 2002.

\bibitem[PV79]{PaVa79}
G.~Papanicolaou and S.R.S. Varadhan.
\newblock Boundary value problems with rapidly oscillating random coefficients.
\newblock In {\em Colloquia Mathematica Societatis J\'{a}nos Bolay}, volume~27,
  1979.

\end{thebibliography}
\end{document}